\documentclass{article}
\usepackage{amsthm, fullpage, amsfonts, url}
\usepackage{amscd,amsmath,amssymb}

\newtheorem{definition}{Definition}
\newtheorem{theorem}{Theorem}
\newtheorem{lemma}{Lemma}
\newtheorem{proposition}{Proposition}
\newtheorem{corollary}{Corollary}

\def\Im{\operatorname{Im}}
\def\de{\delta}
\def\De{\Delta}

\def\si{\sigma}
\def\Si{\Sigma}
\def\al{\alpha}
\def\Ga{\Gamma}
\def\ga{\gamma}
\def\La{\Lambda}
\def\la{\lambda}
\def\Om{\Omega}
\def\om{\omega}
\def\ti{\tilde}
\def\kappa{\varkappa}

\def\Bbb{\mathbb}
\def\SL{\operatorname{SL}}
\def\EXP{\operatorname{EXP}}
\def\LOG{\operatorname{LOG}}
\def\E{{\mathbb E}}
\def\sn{{\mathfrak S}_n}
\def\Prob{\operatorname{Prob}}
\def\inv{\operatorname{Inv}}
\def\supp{\operatorname{supp}}
\def\AUT{\operatorname{AUT}}
\def\SAUT{\operatorname{SAUT}}
\def\Aut{\operatorname{Aut}}

\def\Id{\operatorname{Id}}
\def\be{\begin{equation}}
\def\ee{\end{equation}}

\urldef\vershik\url{vershik@pdmi.ras.ru}
\urldef\natalia\url{natalia@pdmi.ras.ru}
\author {A.~M.~Vershik\thanks{%
St.~Petersburg Department of Steklov Institute of Mathematics.
E-mail: \vershik, \natalia.
Supported in part by RFBR, grant
02-01-00093, by the Russian President grant 2251.2003.1
for Support of the Leading Scientific Schools,
and by the Netherlands Organisation for 
Scientific Research (NWO).%
} \and N.~V.~Tsilevich\footnotemark[1]}
\title {Fock factorizations, and decompositions of the
$L^2$ spaces over general L\'evy processes}
\date{}

\begin{document}
\maketitle

\abstract{
We explicitly construct and study an isometry between the spaces
of square integrable functionals of an arbitrary L\'evy process and
a vector-valued Gaussian white noise. In particular, we obtain explicit
formulas for this isometry at the level of
multiplicative functionals and at the level of orthogonal
decompositions, as well as find its kernel. We consider in detail the
central special case: the isometry between the $L^2$ spaces over 
a Poisson process and the corresponding white noise.
The key role in our considerations is
played by the notion of measure and Hilbert factorizations
and related notions of multiplicative and additive functionals
and logarithm. The obtained
results allow us to introduce a canonical Fock structure
(an analogue of the Wiener--Ito decomposition) in the $L^2$ space
over an arbitrary L\'evy process. An application to the representation
theory of current groups is considered. An example of a non-Fock
factorization is given.
}

\tableofcontents

\section{Introduction: setting of the problem and main results}
\label{sect:intr}

\subsection{Subject of the paper}

This paper is a survey and an exposition of new results
in the field, which has been for a long time relating the classical probability theory,
functional and classical analysis, and combinatorics,
as well as some areas of theoretical physics (the second
quantization, Fock space). We mean the theory of random processes with independent
values (or, in more traditional probabilistic setting,  with  
independent increments) and decompositions of functional spaces
over these processes. 
Such processes may be regarded as a {\it continual generalization
of the notion of a sequence of independent random variables}.
The theory of these processes is closely related to the theory of
infinitely divisible distributions on the real line. It passed a long
way from the original pioneering works by B.~de Finetti and A.~N.~Kolmogorov,
who suggested a formula for infinitely divisible distributions on the line
with a finite variance, subsequent papers of the middle 30s by P.~L\'evy
and A.~Ya.~Khintchin, who proved a general formula for these distributions,
up to the notion of generalized random
processes in the sense of Gelfand--It\^o, which provided a solid base 
for understanding what is a process with independent values.

The central example is of course the Wiener process
(or the Gaussian white noise if we consider generalized processes with
independent values). The measure in the space   $C([0,1])$ of realizations
of this process was described by N.~Wiener in the early 20s;
in the sequel, this measure remained at the center of the whole
stochastic analysis and the theory of stochastic differential equations.
This theory was started already in the 50s in the works by K.~It\^o, which were
continued by many other mathematicians, and by now it has
a vast range of applications in the theory of random processes and 
other fields.

But there is another aspect, which we will consider below and which
is directly related to another source of the theory of ``continual products of
independent random variables'', which is less evident but perhaps the most
important; we mean mathematical physics.
It is worth recalling that the Wiener process is obviously a mathematical
version of the Brownian motion and Einstein--Smolukhovsky process.
However, it was not until the 50s that another remarkable fact became
clear: the so-called Wiener--It\^o--Cameron--Martin
orthogonal decompositions in the Hilbert space of square integrable
functionals of the Wiener processes, whose theory was constructed in
\cite{Wien38, Ito51, CamMart47}, is nothing else but a reproduction of the 
second quantization scheme, which was first suggested by V.~A.~Fock 
in the early 30s and developed in dozens of mathematical and semimathematical
papers of his colleagues; the so-called Fock space, which serves as a base for
constructions of the quantum field theory, representation theory of many
(especially infinite-dimensional) groups, many algebraic constructions, etc.,
has the ``exponential'' structure and orthogonal decompositions into
multiparticle subspaces, exactly as the $L^2$ space 
over the white noise has the decomposition into ``chaoses''
of different orders; and the so-called Wick regularization is
merely the process of orthogonalization of polynomial functionals
of Hermite type (see, e.g., the monograph \cite{Simon76} or the survey
   \cite{DobrMinl77}). Among numerous books related to the subject
under consideration, we mention
\cite{Par, DeMe, Ner}. An important role in the popularization of the Fock
space among the Soviet mathematicians was played by 
F.~A.~Berezin's works (especially on the fermion Fock space), see, e.g.,
\cite{Ber}.

\subsection{Structure of a factorization}   
   
A more careful analysis shows that both Hilbert spaces (the Fock space
and the $L^2$ space over the Wiener measure) {\it have a structure of}
the so-called {\it ``factorization''}, or a continuous tensor
product. It is this structure that corresponds to the intuitive notion of the
``continual product of independent variables''; the existence of such structure
means that the Hilbert
space and the algebra of operators in this space admit infinitely
divisible decompositions into tensor products. This structure of 
decompositions into tensor factors, or factorization, appears not only in the
probability theory, but also in the representation theory of current
groups and fields of $C^*$-algebras, models of the field theory, algebra, etc. 
It goes back to the pre-war works by von Neumann on tensor products
\cite{MurNeum36, Neum38} and was investigated in the paper by Araki and Woods
in the 60s \cite{AW66}.  Its metric (probability-theoretical) counterpart
is more recent, it was suggested by Feldman
\cite{Feld71}. Below we give the definition of a {\it measure factorization}
and a short survey of few papers where it was considered.
Roughly speaking, a continuous measure factorization of a measure space
is a coherent family of decompositions of this space into the direct
product of arbitrarily many measure spaces. The space
of realizations of each process with independent values has
such structure.

{\it The main result of this work says that the Hilbert space of square
integrable functionals over a random process with independent values
in an arbitrary vector space (in short, L\'evy process, though this term is not
quite correct)
has the structure of a Fock factorization. Thus, from this viewpoint,
an arbitrary  L\'evy process has the same
factorization structure 
as the Gaussian process of an appropriate dimension. This dimension is
the only invariant of the factorization up to isomorphism, and
it depends only on the number of points in the support of the L\'evy measure,
which implies that two L\'evy processes with the same cardinality of the
supports of the L\'evy measures determine isomorphic factorizations.}
We present explicit formulas for the factorization-preserving 
isometry between the corresponding $L^2$ spaces. In particular, for any
L\'evy process, one can obtain an orthogonal decomposition, in the
space of square integrable functionals, similar to the classical 
Wiener--It\^o decomposition into ``chaoses''.

Though there are numerous (mostly technical) papers devoted to the construction of
stochastic integrals and analogues of the Wiener--It\^o decomposition
for L\'evy processes (see, e.g., 
\cite{Og72, SegKail76, Eng82}, and also
\cite{Kab75, KalSz89, Sz91, NuSch00, McK73, IsSa83, KSSU98, KL00, BM01}),
the crucial observation that the corresponding
factorizations are isomorphic
remained up to now in the background. This fact is closely related
to the remarkable Araki--Woods theorem \cite{AW66} (though does not follow
from it), which claims that a
factorization having sufficiently many multiplicative vectors
is isomorphic to a Fock factorization, see \cite{AW66, TsirVer98}.
In probabilistic terms, this condition means that the 
above-mentioned existence of the canonical
Fock structure (or Wiener--It\^o decomposition) in the $L^2$ space over an
arbitrary L\'evy process follows from the totality of the set of
multiplicative functionals of the process\footnote{Recall that 
a subset of a Hilbert
space is called total if its linear span is everywhere dense, 
and that a functional
of a process defined on a set $X$
is called multiplicative if its value at an arbitrary 
realization of the process
is equal to the product of its values at the restrictions of this 
realization to subsets forming an arbitrary finite measurable 
partition of the whole
set $X$ (see Sect.~\ref{sect:fact} for a precise definition).}.
A multiplicative functional of
a Gaussian process is the exponential of a linear functional,
however, for general L\'evy processes, the set of multiplicative
functionals is much wider (see Sect.~\ref{sect:log}).
The established Fock structure in the $L^2$ space over an arbitrary
L\'evy process and the isometry between these spaces and the $L^2$ spaces over
Gaussian processes make unnecessary numerous special constructions of 
orthogonal decompositions (and stochastic integrals) in each particular case.
Another important detail is that the base space
over which the process is defined is irrelevant for these issues,
it may
be an arbitrary measure space rather than an interval or the line 
as usual. In particular, the isometry under consideration
applies to random fields. The only
advantage of the one-dimensional situation is that 
in this case one may argue in terms of 
processes with independent increments rather than independent values,
and consider the Wiener process instead of the white noise, which is 
sometimes more convenient. However, we consider L\'evy processes over
an arbitrary base space.

\subsection{Isomorphism of factorizations. Logarithmic operation
in factorizations. Kernel}

The simplest example of the isometry under
consideration is the isometry between the $L^2$ spaces over the Poisson and Wiener 
processes. An analogy between the orthogonal structures in these spaces
was observed in many papers; however, the existence of an
isometry was not established even in this case. It is worth mentioning 
that the existence of this
isometry was originally obtained in \cite{GGV75}
from the equivalence of two
realizations of the canonical representation of the groups of diffeomorphisms; 
later, this isometry was studied in \cite{Ner97} in terms of the so-called
holomorphic model of the Fock space (i.e., the Fock space
realized as the Hilbert space of
holomorphic functionals rather than the $L^2$ space), which is more popular among
physicists. Final explicit formulas for the isometry between
the $L^2$ spaces over the Poisson and Gaussian processes are apparently new.

The general case of the isometry between an arbitrary L\'evy process and a
Gaussian process of an appropriate dimension can be
reduced to the above-mentioned Poisson--Gauss case by means of
the L\'evy--Khintchin decomposition or Poissonian construction
of L\'evy processes. But in this case we need to consider
vector-valued white noises.
Namely, it is natural to take the Hilbert space $L^2(\mathbb R, \Pi)$, 
where $\Pi$  is the L\'evy measure of the L\'evy process under consideration,
as the space of values of the white noise.
In particular, in the case of a 
Poisson process this space is one-dimensional, hence the Poissonian
factorization is isomorphic to the factorization generated by the
ordinary (one-dimensional) Wiener process.

As mentioned above, the only Hilbert invariant
(i.e., invariant up to arbitrary isometries of the Hilbert space)
of the factorizations arising in the theory
of L\'evy processes  is the dimension of the 
space $L^2(\Pi)$, i.e., the number of points in the support of the L\'evy 
measure. For example, the stable processes and the gamma processes,
which are of importance for applications,
generate the same factorization as the Gaussian white noise with values in
an infinite-dimensional Hilbert space. At the same time, the metric invariants
of factorizations, 
i.e., invariants up to measure-preserving transformations, are much 
more detailed; namely, as shown by Feldman \cite{Feld71},
two measure factorizations 
generated by L\'evy processes with L\'evy measures $\Pi_1$ and $\Pi_2$  
are isomorphic if and only if  $(\mathbb R, \Pi_1)$ and $(\mathbb R, \Pi_2)$ 
are isomorphic as measure spaces.
The  metric classification of processes plays an important role
in the theory of decreasing families of $\si$-fields (filtrations)
(see \cite{V94}).

The explicit construction of the Fock--Wiener--It\^o structure in the
space of square integrable functionals 
over an arbitrary L\'evy process, i.e., the decomposition
of this space into the orthogonal sum of the symmetric tensor powers of the first
chaos, is based on a kind of 
{\it ``taking logarithm'' of multiplicative functionals},
resulting in the space of additive 
functionals, i.e., the first chaos. 
However, the rule for calculating this ``logarithm''  
substantially depends on the factorization, and in general it does not coincide
with taking
the ordinary logarithm 
(the coincidence takes place only for Gaussian
processes; in this case the set of additive functionals coincides
with the set of linear functionals of the process). The existence of
this ``logarithm'' for general L\'evy processes is not obvious;
the proof of the Araki--Woods theorem
consists essentially in constructing
this logarithm (in a slightly more general context); a more
explicit version of this construction can be found in \cite[Appendix A]{TsirVer98}.
Below (Sect.~\ref{sect:log}) we calculate the logarithm
for the Gaussian and Poisson cases (which,
as we have mentioned above, exhaust the general case of L\'evy processes).

In order to determine uniquely the isometry, it suffices to establish a correspondence
between the sets of multiplicative functionals or 
between the linear subspaces of additive functionals (``first
chaoses''); if we fix a bijection between the canonical bases in these
subspaces, then the isometry is unique. {\it The isometry can be also defined by 
a kernel}, i.e., a generalized function in realizations of the two processes.
This kernel for the Poisson--Gauss case is computed in 
Sect.~\ref{sect:kernel}; remarkably, 
it is defined in purely combinatorial terms.
The correspondence of successive chaoses in this case reduces to
a correspondence between the Hermite and Charlier functionals.
The kernel is not positive, hence the isometry is neither
Markovian, nor multiplicative operator.

Let us say some words on the combinatorial and
analytical aspects of the problem. The isometry under consideration,
for a L\'evy process with L\'evy measure $\Pi$, takes a
very explicit form
if the moment problem for the measure $t^2d\Pi(t)$ is definite.
In this case we can consider the basis of orthogonal polynomials in 
$L^2(\mathbb R,t^2d\Pi(t))$ and obtain interesting relations between orthogonal 
polynomials with respect to various measures on the real line
(see, e.g., (\ref{hermcharid})). The most important role is played
by the Hermite and Charlier polynomials.
Combinatorial aspects of constructing stochastic integrals
for a wide class of processes were most explicitly considered in \cite{RotaW97}.

\subsection{Relation to the representation theory}

One of the most important applications of the 
isometries between Hilbert spaces
of functionals over various L\'evy processes is the representation theory
of infinite-dimensional groups, namely, of current and gauge groups, and of
Kac--Moody algebras. Unitary representations of these groups naturally
generate Hilbert factorizations. Usually, these representations are
realized in the Fock space.
The model of the Fock space 
as the $L^2$ space over the white noise or the close holomorphic model
of this space as the space of holomorphic
functions in infinitely many variables are only two of possible models.
If we fix an arbitrary commutative subgroup (subalgebra) of the 
infinite-dimensional group (algebra), and construct the representation
of this algebra
where this subgroup (subalgebra) is diagonalized,
then we obtain immediately one of these models. For instance,
the commutative model of the
canonical representation of the group of 
$\operatorname{SL}(2,\mathbb R)$-currents 
(see \cite{GGV73})
with respect to  the subgroup of unipotent matrices yields
the isometry of the $L^2$ spaces over the
infinite-dimensional white noise and the gamma
process (see \cite{GGV85}).
In particular, this observation led to discovering new
symmetry properties of the gamma process.
This example is considered in detail in \cite{TVY01}.

\subsection{Further development}

Processes with independent values and factorizations appear in 
much wider context than discussed above. First of all, one may 
consider an analogue of Wiener and other processes on manifolds,
groups, semigroups, and even more general systems; the only thing
we need is a distribution that is infinitely divisible with
respect to a composition of measures. Moreover, the notion of a
composition of measures may be very general and even not related to a
group or semigroup law.

For example, consider the simplest nonlinear case --- the Brownian
motion on the sphere, or the rotation group
$SU(n)$. It was proved already in the 60s (see \cite{McK60})
that this process is linearizable, i.e., it can be represented by means
of the ordinary Wiener process; thus the corresponding factorization
is a Fock factorization. The method of stochastic differential equations,
used to obtain this result, does not give an answer to the question
investigated in this paper: what is the decomposition into ``chaoses''
in this case, i.e., what is the explicit isomorphism of the corresponding
spaces of functions. Note that an answer to this question would provide
a direct proof of the linearizability of the Brownian motion.
The strongest result on linearizability was recently obtained by
B.~Tsirelson \cite{Tsir98}. It claims that a weakly continuous Brownian
motion on the unitary group of an infinite-dimensional Hilbert space
is linearizable; it follows that the same holds for all groups that have
a faithful unitary representation. The question of finding an explicit isometry
and an explicit decomposition of the $L^2$ space into chaoses remains open.
But for a wider class of groups, for example, for the group of isometries
of the universal Urysohn space or the group of homeomorphisms of a compact
space, even the linearizability of the Brownian motion is still not
proved. And it is absolutely unclear if this is true for

a) non-Gaussian processes with independent group values, for example,
L\'evy processes on finite-dimensional or infinite-dimensional
groups.

b) any processes with independent semigroup and more general values.

Of great interest are questions concerning the {\it metric} classification
of the factorizations generated by L\'evy processes with arbitrary
(nonlinear) values; the simplest of these questions is whether the factorization generated
by the Brownian motion on the two-dimensional sphere is {\it metrically}
isomorphic to the one-dimensional Gaussian factorization.

\subsection{Non-Fock factorizations}
\label{sect:nonfock}

The new stage of the development of the theory of factorizations
is related to deeper questions.

In \cite{Feld71}, the following question, which goes back to 
S.~Kakutani, was discussed: whether it is true that every measure factorization is isomorphic
to a Fock factorization, or, in another terminology, is linearizable?
We will discuss this question in more detail at the end of the paper; 
here we only  mention that, in view of a theorem similar to 
the remarkable Araki--Woods theorem on Hilbert factorizations \cite{AW66} 
(see also Sect.~\ref{sect:fockspace} below),  
this question is equivalent to the question if
there are sufficiently many factorizable (multiplicative) functionals
in this measure space.
It turns out that the cases are possible when 
there are no multiplicative functionals
except constants, and such examples of non-Fock factorizations (``black
noise'') with characteristic strong nonlinearity were constructed in
\cite{TsirVer98} for a base of dimension
0 and 1. The constructed examples are in no sense generalized random
processes with group or semigroup values. We give a short version of the
example for a base of dimension 0 in Appendix A.

There is another important difference of 
these factorizations from the Fock ones: unlike Fock factorizations, which
are defined on the complete Boolean algebra
of classes of $\bmod0$ coinciding measurable sets, these factorizations
are defined on a more narrow Boolean algebra; it was this fact that
caused the restriction
on the dimension of the base. As shown by Tsirelson \cite[Sect.~6c]{Ts03}, 
on the complete Boolean algebra, every 
factorization satisfying certain continuity conditions is a Fock factorization.

The problem of applying non-Fock factorizations in the representation
theory of current groups, fields, and $C^*$-algebras, and in the quantum theory
is still actual; the existence of such factorizations apparently opens
new possibilities in the representation theory of infinite-dimensional and
field objects.

\subsection{Structure of the paper}

The paper is organized as follows.

\S\ref{sect:basicdef} contains the necessary background on 
factorizations and processes with
independent values.
We also give the  definition
of the logarithm determined by a factorization. The main results are contained in
\S\ref{sect:isom0}, where we consider in detail the fundamental special case ---
the canonical isometry between the spaces of square integrable 
functionals over the Poisson and 
Gaussian processes. For both processes, we compute all the 
above-mentioned characteristics
(multiplicative and additive functionals, logarithm, 
orthogonal decomposition into stochastic integrals). As model cases,
we consider the finite-dimensional analogues corresponding to
a finite base space.
We would like to draw the reader's attention to 
formula (\ref{hermcharid})
for the classical Hermite and Charlier orthogonal polynomials,
which is apparently new. The problem of determining the kernel
of the canonical isomorphism
(formula (\ref{kernel})) is new both in setting and in suggested solution.
There are several proofs of this formula
including a purely combinatorial one. 
In \S\ref{sect:genisom}, we study the isometry for
a general L\'evy process; this requires no
substantially new ideas, since a well-known construction allows one to
represent such process as a Poisson process on a wider space. Though this
representation is well-known, nevertheless it was
not realized that in this general case
the $L^2$ space over an arbitrary
L\'evy process with L\'evy measure $\Pi$
is isometric to the $L^2$ space over a $L^2(\mathbb R,\Pi)$-{\it valued} 
Wiener process. 
Thus all formulas in the general case merely
reproduce the corresponding formulas for the Poisson--Gauss case. 
Finally, in \S\ref{sect:repr}, we
consider an example of 
applying the isometry between the Fock space and
the space of square integrable
functionals over the gamma process. In fact, it was this example that
gave rise to the series of papers 
\cite{GGV73, GGV74, GGV75, GGV83, GGV85, TVY01} resulting
in the present understanding of the whole situation.

\bigskip

The authors are grateful to Professor M.~Yor for a number of important 
bibliography references concerning the theory of L\'evy processes of 
former years.

\section{Basic definitions}
\label{sect:basicdef}
\subsection{Factorizations}
\label{sect:fact}

We will consider two types of factorizations: {\it Hilbert factorizations}
and {\it measure factorizations}.

Given a Hilbert space $\cal H$, denote by $\cal B(\cal H)$ the algebra of
all bounded linear operators on $\cal H$, and let
${\cal R}(\cal H)$ be the lattice
of all von Neumann algebras on $\cal H$ (recall that the lattice operations
in ${\cal R}(\cal H)$ are defined as follows:
$R_1\wedge R_2=R_1\cap R_2$ and $R_1\vee R_2=(R_1\cup R_2)''$, where
$R'=\{a\in{\cal B}({\cal H}):\,ar=ra\;\forall r\in R\}$ is the commutant of 
$R$; for an exposition of the theory of von Neumann algebras, see, e.g., 
\cite{Naimark})). 
The following definition of a Hilbert factorization goes back to von Neumann
\cite{MurNeum36, Neum38}.

\begin{definition}
\label{def:hilbfact}
A (type I) {\em Hilbert factorization} of a Hilbert space
$\cal H$ over a Boolean algebra $\cal A$
is a map $\xi:{\cal A}\to{\cal R}(\cal H)$ such that each 
algebra of operators $\xi(A)$ is a type I factor, and for all 
$A, A_1,A_2,\ldots\in\cal A$, the following conditions hold:
\begin{itemize}
\item $\xi(A_1\wedge A_2)=\xi(A_1)\wedge\xi(A_2)$;
\item $\xi(A_1\vee A_2)=\xi(A_1)\vee\xi(A_2)$;
\item $\xi(A')=\xi(A)'$;
\item $\xi(0_{\cal A})=\{\alpha\cdot\Id_{\cal H},\,\alpha\in\mathbb C\}=
1_{\cal H}$\footnote{The algebra $\{\alpha\cdot\Id_{\cal H},\,\alpha\in\mathbb C\}$
is traditionally denoted by $1_{\cal H}$, though
it is the {\em zero} of the lattice of von Neumann algebras in $\cal H$.},
where $\Id_{\cal H}$ is the identity operator in $\cal H$;
\item $\xi(1_{\cal A})={\cal B}(\cal H)$
\end{itemize}
A {\em factorized Hilbert space} $({\cal H}, \xi)$
is a Hilbert space $\cal H$ equipped with a Hilbert factorization $\xi$
of its operator algebra. The Boolean algebra $\cal A$ is called
the {\em base} of the factorization.
\end{definition}

\noindent{\bf Remark.} If $({\cal H}, \xi)$ is a factorized Hilbert space,
then, as shown in \cite{AW66}, for each $A\in\cal A$,
there is a subspace
${\cal H}_A\subset{\cal H}$ such that for each finite partition
$A_1,\ldots,A_n$ of the unity element $1_{\cal A}$ of the Boolean 
algebra ${\cal A}$, we have
${\cal H}={\cal H}_{A_1}\otimes\ldots\otimes{\cal H}_{A_n}$ and
$\xi(A_k)=1_{{\cal H}_{A_1}}\otimes\ldots\otimes
1_{{\cal H}_{A_{k-1}}}\otimes{\cal B}({\cal H}_{A_k})
\otimes 1_{{\cal H}_{A_{k+1}}}\otimes
\ldots\otimes1_{{\cal H}_{A_n}}$.
\smallskip

The notion of a measure factorization was introduced by Feldman \cite{Feld71}.
We follow the presentation adopted by Tsirelson and Vershik \cite{TsirVer98}.

\begin{definition}
\label{def:measfact}
Let $(\Om,{\frak A},\mathbb P)$ be a probability space
(which is always assumed to be a continuous Lebesgue space).
Denote by $\Si(\mathbb P)$ the complete lattice
of all sub-$\si$-fields (containing all negligible sets) of the $\si$-field
$\frak A$.
A {\em measure factorization} of $(\Om,{\frak A},\mathbb P)$
over a Boolean algebra $\cal A$ is a map 
$\zeta:{\cal A}\to\Si(\mathbb P)$ 
such that for all $A,A_1,A_2,\ldots\in\cal A$,
the following conditions\footnote{These conditions are not independent.} hold:
\begin{itemize}
\item $\zeta(A_1\wedge A_2)=\zeta(A_1)\wedge\zeta(A_2)$;
\item $\zeta(A_1\vee A_2)=\zeta(A_1)\vee\zeta(A_2)$;
\item $\zeta(A')$ is an independent complement\footnote{Such complement
is not unique.} of the $\si$-field
$\zeta(A)$, i.e., $\zeta(A)\wedge\zeta(A')=0$, $\zeta(A)\vee\zeta(A')=1$,
and the  $\si$-fields $\zeta(A)$ and $\zeta(A')$ are independent
(which means that  $\mathbb P(E_1\cap E_2)=\mathbb P(E_1)\mathbb P(E_2)$
for all $E_1\in\zeta(A)$ and $E_2\in\zeta(A')$);
\item ${\zeta}(0_{\cal A})={\frak A}_0$ (the trivial $\si$-field);
\item $\zeta(1_{\cal A})=\frak A$.
\end{itemize}
A {\em factorized measure space} $(\Om,{\frak A},\mathbb P,\zeta)$
is a probability space equipped with
a measure factorization $\zeta$ over some Boolean algebra $\cal A$. 
The Boolean algebra $\cal A$ is called
the {\em base} of the factorization.
\end{definition}

In this paper, we will consider only {\it continuous} 
(Hilbert and measure) factorizations in the sense of the following
condition (which is called ``minimal up continuity condition'' in
\cite{TsirVer98}). In what follows, the term ``factorization''
means ``continuous factorization'', unless otherwise stated.

\begin{definition}
A Hilbert factorization $\xi$ (respectively, a measure factorization
$\zeta$) over a Boolean algebra
$\cal A$ is called {\em continuous} if
$\bigvee_{A\in S}\xi(A)={\cal B}(\cal H)$
(respectively, $\bigvee_{A\in S}\zeta(A)=\frak A$)
for every maximal ideal $S\subset\cal A$.
\end{definition}

The most important examples are factorizations over the Boolean algebra
$\frak B$ of {\it all Borel sets} of a standard Borel space $(X,\frak B)$
(which will be called {\it factorizations over the Borel space $(X,\frak B)$}), and
factorizations
over the Boolean algebra of {\it $\bmod0$ classes of measurable subsets} of 
a Lebesgue space $(X,\nu)$ (which will be called {\it factorizations over
the Lebesgue space $(X,\nu)$}).
In this paper, we consider only factorizations of these two types.
Moreover, in all our examples, a factorization over a Borel space
can be correctly extended to a factorization over the corresponding Lebesgue space
(in fact, this is a consequence of the fact that all factorizations considered
in this paper turn out to be Fock factorizations, see below;
in the case of non-Fock factorizations (i.e., of factorizations
that are not isomorphic to Fock ones),
the base Boolean algebra is more narrow than the algebra of all
Borel sets, see Appendix A).

\begin{definition}\label{def:isomfact}
1) Two {\em factorized Hilbert spaces} $({\cal H}_1,\xi_1)$ and $({\cal H}_2,\xi_2)$
over Boolean algebras ${\cal A}_1$
and ${\cal A}_2$, respectively, are called {\em isomorphic}
if there exists an isomorphism of Boolean algebras
$S:{\cal A}_1\to{\cal A}_2$ and
an isometry of the Hilbert spaces $T:{\cal H}_1\to{\cal H}_2$ such
that the following diagram is commutative:
$$
\begin{CD}
{\cal A}_1@>S>>{\cal A}_2\\
@VV\xi_1V     @VV\xi_2V\\
{\cal R}({\cal H}_1)@>\bar T>>{\cal R}({\cal H}_2).
\end{CD}
$$
Here $\bar T$ is the operator from ${\cal R}({\cal H}_1)$
to ${\cal R}({\cal H}_2)$ generated by the isometry $T$ of the Hilbert spaces.

2) In a similar way,
two {\em factorized measure
spaces} $(\Om_1,{\frak A}_1,\mathbb P_1,\zeta_1)$ and 
$(\Om_2,{\frak A}_2,\mathbb P_2,\zeta_2)$ over Boolean algebras ${\cal A}_1$
and ${\cal A}_2$, respectively, are called {\em 
isomorphic} if there exists an isomorphism of Boolean algebras
$S:{\cal A}_1\to{\cal A}_2$ and 
an isomorphism of measure spaces 
$T:(\Om_1,{\frak A}_1,\mathbb P_1)\to (\Om_2,{\frak A}_2,\mathbb P_2)$
such that the following diagram is commutative:
$$
\begin{CD}
{\cal A}_1@>S>>{\cal A}_2\\
@VV\zeta_1V     @VV\zeta_2V\\
\Si(\mathbb P_1)@>T>>\Si(\mathbb P_2).
\end{CD}
$$
\end{definition}

\begin{definition}
If $T$ is an isomorphism of (Hilbert or measure)
factorizations defined over the same
Boolean algebra $\cal A$, and the corresponding automorphism $S$ is the
identity automorphism of the Boolean algebra $\cal A$
(i.e., the base is fixed), 
then $T$ is called a {\em special} isomorphism of factorizations.
\end{definition}

The following lemma is obvious.

\begin{lemma} \label{l:meas2hilb}
Each measure factorization $(\Om,{\frak A},\mathbb P,\zeta)$
over a Boolean algebra $\cal A$ generates  
a Hilbert factorization in the space ${\cal H}=L^2(\Om,\mathbb P)$
with the same base.
Namely, for each $A\in\cal A$, let
$H_A=L^2(\Om,\zeta(A),\mathbb P|_{\zeta(A)})\subset L^2(\Om,\frak A,\mathbb P)$.
Then the map $\xi:{\cal A}\to{\cal R}(\cal H)$ given by
\be\label{meas2hilb}
\xi(A)={\cal B}(H_A)\otimes 1_{H_{A'}}
\ee 
is a Hilbert factorization in $\cal H$.
\end{lemma}

In this paper, we deal only with Hilbert factorizations of this type.
Note that nonisomorphic measure 
factorizations may generate isomorphic Hilbert factorizations,
since not every isometry of the $L^2$ spaces is generated by some
isomorphism of the underlying measure spaces.

\smallskip\noindent{\bf Remark.} In fact, a measure factorization is a triple
$(\xi,{\cal Z}, \psi)$, where $\xi:{\cal A}\to{\cal B}({\cal H})$ is
a Hilbert factorization, ${\cal Z}\subset {\cal B}({\cal H})$ is a maximal
commutative subalgebra, and  $\psi\in\cal H$ is a factorizable vector
(see Definition~\ref{def:factvect} below)
such that ${\cal Z}\cap\xi(A)$ is a maximal commutative subalgebra of 
$\xi(A)$ for each $A\in\cal A$, and $\psi$ is $\cal Z$-cyclic.
\smallskip

The key role in the study of factorizations is played by the 
notion of multiplicative and additive functionals.

\begin{definition}\label{def:factvect}
Let $(\Om,{\frak A},\mathbb P,\zeta)$ be a factorized measure space over
a Boolean algebra $\cal A$. A measurable function 
$F:\Om\to\mathbb C$
is called an {\em additive} (respectively, {\em multiplicative}) 
{\em functional} if for every finite 
partition $A_1,\ldots,A_n$ of the unity element $1_{\cal A}$ 
of the Boolean algebra $\cal A$, there exist functions
$F_{A_1},\ldots,F_{A_n}:\Om\to\mathbb C$ such that $F_{A_k}$ 
is $\zeta(A_k)$-measurable
and $F=F_{A_1}+\ldots+F_{A_k}$ (respectively, $F=F_{A_1}\cdot\ldots\cdot F_{A_k}$).
\end{definition}

\begin{definition}
Let $({\cal H},\xi)$ be a Hilbert factorization over 
a Boolean algebra $\cal A$.  A vector
$h\in\cal H$ is called {\em factorizable} 
if for very finite partition $A_1,\ldots,A_n$ of the unity element
$1_{\cal A}$, there exist operators
$P_k\in\xi(A_k)$ such that the one-dimensional projection
$P_h$ to the vector $h$ can be represented in the form 
$P_h=P_1\otimes\ldots\otimes P_n$. Alternatively, there exist vectors
$h_{A_i}\in H_{A_i}$ (where $H_{A_i}$ are the subspaces from the remark after
Definition~\ref{def:hilbfact}) such that
$h=h_{A_1}\otimes\ldots\otimes h_{A_n}$.

Analogously, a vector $h\in\cal H$ is called {\em additive} if for any
finite partition $A_1,\ldots,A_n$ of the unity element $1_{\cal A}$,
there exist vectors $h_{A_i}\in H_{A_i}$ such that
$h=h_{A_1}+\ldots+h_{A_n}$.
\end{definition}

If a factorized Hilbert space $(L^2(\Om, \mathbb P),\xi)$ 
is generated by a factorized measure space
$(\Om,{\frak A},\mathbb P,\zeta)$ as in Lemma~\ref{l:meas2hilb}, 
then the set of factorizable vectors in $L^2(\Om, \mathbb P)$
coincides with the set of
square integrable multiplicative functionals in $(\Om,{\frak A},\mathbb P,\zeta)$.
Since in this paper we consider only Hilbert factorizations of this type, we will
use the term ``multiplicative functionals'' for factorizable vectors in
$L^2(\Om,\mathbb P)$.
The set of additive vectors in a factorized Hilbert space 
${\cal H}$ is a linear subspace (maybe zero), and in the case
${\cal H}=(L^2(\Om, \mathbb P),\xi)$ it coincides with the set of
square integrable additive functionals in 
$(\Om,{\frak A},\mathbb P,\zeta)$.

\subsection{First examples: Fock factorizations, Gaussian and Poisson processes}

\subsubsection{Fock spaces and Fock factorizations}
\label{sect:fockspace}

The (boson) {\em Fock space} $\EXP H$ over a Hilbert space $H$ is the
symmetrized tensor exponential
$$
\label{Fock}
\EXP H=S^0H\oplus S^1H\oplus\ldots\oplus S^nH\oplus\ldots,
$$
where $S^nH$ is the $n$th symmetric tensor power of $H$.  Given $h\in H$, let
$$
\EXP h=1\oplus h\oplus\frac{1}{\sqrt{2!}}h\otimes h\oplus\frac{1}{\sqrt{3!}}h\otimes h\otimes h\oplus\ldots.
$$
The vectors $\{\EXP h\}_{h\in H}$ are linearly independent, and
$$
(\EXP h_1,\EXP h_2)_{\EXP H}=\exp(h_1,h_2)_H,
$$
where $(\cdot,\cdot)_{\EXP H}$ stands for the scalar product in $\EXP H$.
In particular,
$$
\|\EXP h\|^2=\exp\|h\|^2.
$$

The principal example of a Hilbert factorization is given by the following
construction.

\begin{definition} Given a Lebesgue space $(X,\nu)$,
consider the direct integral of Hilbert spaces
${\cal K}=\int^{\oplus}K(x)d\nu(x)$ and the corresponding Fock space 
${\cal H}=\EXP \cal K$. For each measurable $A\subset X$, let
$H_A=\EXP {\cal K}(A)$, where ${\cal K}(A)=\int^{\oplus}_AK(x)d\nu(x)$, and
set $\xi(A)={\cal B}(H_A)\otimes 1_{H_{A'}}$.
The obtained Hilbert factorization $(\cal H,\xi)$ is called a
{\em Fock factorization}.
\end{definition}

In particular, if $\dim K(x)\equiv 1$, then
${\cal H}=\EXP L^2(X,\nu)$, and
$H_A=\EXP L^2(A,\nu_A)$, where $\nu_A$ is the restriction of the measure 
$\nu$ to the subset
$A$. More generally, if $\dim K(x)\equiv n$ ($n=1,2,\ldots,\infty$), 
then $\cal H$ can be identified with $\EXP L^2((X,\nu);H)$, 
where $H$ is a Hilbert space of dimension $n$, and $L^2((X,\nu);H)$ is the space
of square integrable $H$-valued functions on $(X,\nu)$. 
The corresponding factorization is called a
{\em homogeneous Fock factorization of dimension $n$}.

The set of multiplicative (factorizable) vectors in a Fock space
${\cal H}=\EXP{\cal K}$ is
$$
{\cal M}=\{c\cdot\EXP h,\;h\in {\cal K},c\in\mathbb C\}.
$$
The linear subspace of additive vectors in the Fock space ${\cal H}=\EXP{\cal K}$
can be identified with the space
${\cal K}=\int^{\oplus}K(x)d\nu(x)$,
where ${\cal K}$ is embedded in ${\cal H}=\EXP{\cal K}$ 
as the subspace of the first chaos:
${\cal K}\ni h\mapsto 0\otimes h\otimes0\otimes\ldots\in\EXP\cal K$.
Thus the space of the first chaos, as well as the set of multiplicative
functionals, is defined in invariant terms; it has the structure of the direct 
integral of Hilbert spaces, the base of the integral coinciding obviously
with the base of the factorization.

Fock factorizations are characterized by the following
important theorem.

\begin{theorem}[Araki--Woods \cite{AW66}] 1) A Hilbert factorization $({\cal H},\xi)$
over a nonatomic Boolean algebra $\cal A$ is a Fock factorization if
and only if the set of factorizable vectors is total in $\cal H$.

2) The complete invariant of a Fock factorization is the set of values
assumed by the dimension function
$\dim K(x)$ at sets of positive measure:
$\{n\in\{0,1,\ldots;\infty\}:\;\nu(\{x:\dim K(x)=n\})>0\}$. Thus
two Fock factorizations are isomorphic if and only if their dimension
functions are equivalent in the following sense: the set of values
that they assume at sets of positive measure coincide.
\end{theorem}

\subsubsection{Gaussian white noise and Gaussian factorizations}
\label{sect:gauss}

Consider the {\em Gaussian white noise} $\al$ on a Lebesgue 
space $(X,\nu)$, i.e., the generalized random 
process\footnote{The notion of a generalized random process was developed
by I.~M.~Gelfand \cite{Gel55} (see also \cite{GelVil61})
and K.~It\^o \cite{Ito54}. Concerning the Gaussian white noise, see also
\cite{Skor} and \cite{Neveu}.}
on the Hilbert space $L^2(X,\nu)$ with the characteristic functional given
by the formula
\be
\E e^{i\langle h,\cdot\rangle}=e^{-\frac12\|h\|^2},
\quad h\in L^2(X,\nu),
\ee
where $\|h\|$ is the norm of a vector $h$ in the space $L^2(X,\nu)$.
(This process can also be defined explicitly in the nuclear extension
$\hat H$ of the Hilbert space $L^2(X,\nu)$ corresponding to the quadratic form
$B(\cdot,\cdot)=(\cdot,\cdot)$. In particular, the space
$L^2(\al)$ of square integrable functionals of the white noise can be identified
with $L^2(\hat H,\mu)$,
where $\mu$ is the standard Gaussian measure in $\hat H$.)

It is well-known (see, e.g., \cite{Simon76})
that the space $L^2(\al)$ of square integrable
functionals of the white noise is canonically isomorphic to
the Fock space $\EXP H$, where $H=L^2(X,\nu)$. This isomorphism
is given by the formula
\be
\label{fockgauss}
\EXP h\leftrightarrow e^{-\frac{\|h\|^2}{2}+\langle h,\cdot\rangle}.
\ee
In particular, the vacuum vector $\EXP 0$ corresponds to the unity function
$1\in L^2(\al)$, and 
the $n$-particle subspace $S^nH$ is identified with the subspace of 
$L^2(\al)$ spanned by the $n$-multiple stochastic integrals, i.e., 
by the generalized Hermite functionals of order $n$, see Sect.~\ref{sect:orth0}.
The structure of a unitary ring in
$L^2(\al)$ was described axiomatically in \cite{Ver62-1}, see also \cite{Ver62-2}.

For each measurable set $A\subset X$, denote by $\zeta_\al(A)$ the
$\si$-field generated by the restriction of the process
$\al$ to $A$. It is easy to check that we obtain a measure factorization over
$(X, \nu)$ called a {\it Gaussian} factorization. According to the general
construction of Lemma \ref{l:meas2hilb}, the white noise $\al$
determines also a Hilbert factorization
$\xi_\al$ in the space $L^2(\al)$.
The following well-known proposition is a consequence of the
isomorphism~(\ref{fockgauss}) between $L^2(\al)$ and $\EXP L^2(X,\nu)$.

\begin{proposition}
The Gaussian white noise on an arbitrary Lebesgue 
space $(X,\nu)$ generates a Fock
factorization.
\end{proposition}

In particular, the set of multiplicative functionals in the space $L^2(\al)$
is
$\{c\cdot e^{\langle h,\cdot\rangle},\;h\in L^2(X,\nu),\,c\in\mathbb C\}$, and
the set of additive functionals is
$\{c\cdot \langle h,\cdot\rangle,\;h\in L^2(X,\nu),\,c\in\mathbb C\}$
(see Sect.~\ref{sect:log}).

\smallskip\noindent{\bf Remark.} One may consider the Gaussian process on
$(X,\nu)$ with an arbitrary variance
$\si^2>0$, i.e., the generalized random process
on $L^2(X,\nu)$ with the characteristic functional
$$
\E e^{i\langle h,\cdot\rangle}d\mu(\cdot)=e^{-\frac{\si^2}2\|h\|^2},
\quad h\in L^2(X,\nu).
$$
Clearly, this case reduces to the standard white noise by means of the map
$\langle h,\cdot\rangle
\mapsto\si^{-1}\langle h,\cdot\rangle$, and the factorization determined
by such process is also a Fock factorization.

\subsubsection{Poisson process and Poissonian factorization}
\label{sect:poisson}

A standard reference on the theory of Poisson processes is the book
\cite{Ki93}.

By definition, a point {\it configuration} in the space $X$ is a (non-ordered)
empty, finite, or countable set of points of $X$ with
positive (integral) multiplicities.
Denote by ${\cal E}={\cal E}(X)$ the set of all point configurations
on $X$.

The {\em point Poisson process} on the space $X$ with mean measure $\nu$
(in short, the Poisson process on $(X,\nu)$)
 is a random configuration $\pi\in{\cal E}$ such that
for each measurable subset $A\subset X$, the random variable $\#\{\pi\cap A\}$
has the Poisson distribution with parameter $\nu(A)$, i.e.,
$\Prob\{\#\{\pi\cap A\}=n\}=\frac{\nu(A)^n}{n!}e^{-\nu(A)}$;
and for any disjoint measurable subsets $A_1,\ldots,A_n\subset X$,
the random variables $\#\{\om\cap A_k\}$, $k=1,{\ldots},n$, are independent.

Note that a Poisson process can be regarded as a random measure
$\tau(A)=\#\{\pi\cap A\}$ on the space $X$.
Denote by $\cal P$ the distribution of the Poisson process in the space
of configurations $\cal E$.

As in the Gaussian case, for each measurable subset
$A\subset X$, denote by $\zeta_\pi(A)$
the $\si$-field generated by the restriction of the Poisson process
$\pi$ to  $A$. We obtain a 
{\it Poissonian} measure factorization over
$(X, \nu)$. According to the general
construction of Lemma \ref{l:meas2hilb}, 
the Poisson process $\pi$ determines also a Hilbert factorization
$\xi_\pi$ in the space $L^2(\pi)$.

\subsection{Logarithm}
\label{sect:loggen}

In this section, we describe the operation of ``logarithm'', introduced in 
\cite{TsirVer98}, which allows 
one to construct the space of additive functionals of a factorized measure
space $(\Om,{\frak A},\mathbb P,\zeta)$, 
given the set of multiplicative functionals. For example, 
if the set of multiplicative functionals is total in $L^2(\Om,\mathbb P)$, 
then, by the Araki--Woods theorem,
the corresponding Hilbert factorization of $L^2(\Om,\mathbb P)$
is a Fock factorization, and the logarithmic operation
allows one to construct the space of the first chaos, and hence 
to recover the whole Fock structure in $L^2(\Om,\mathbb P)$
by means of the standard orthogonalization process
(see Sect.~\ref{sect:orthdec}).
The logarithmic operation
depends substantially on the factorization. In general, it does not 
coincide with the ordinary logarithm of a multiplicative functional.

Let $(\Om,{\frak A},\mathbb P,\zeta)$ be 
a factorized measure space over a Lebesgue space $(X,\nu)$.
Denote by $\cal A$ the space of 
all square integrable additive functionals, and by $\cal M$ the 
space of all square integrable multiplicative functionals
on this factorized space. Given $F\in\cal M$,
denote by $F_A(\cdot)$ the $\zeta(A)$-measurable function from the definition
of a multiplicative functional (it is defined uniquely up to scalar factor).
For each measurable subset $A\subset X$, let
$m_F(A)=\log\E|F_A(\cdot)|^2-\log|\E F_A(\cdot)|^2$.
Then $m_F$ is a nonatomic measure on $X$ (\cite[Lemma A2]{TsirVer98}).

\begin{theorem}[\cite{TsirVer98}, Theorem A6] 
\label{th:log}
There is a natural one-to-one
correspondence $\LOG_\zeta:\cal M\to\cal A$. Given $F\in\cal M$ with
$\E F=1$,
\be
\label{log}
\LOG_\zeta F(\cdot)=\lim_{\max m_F(A_k)\to0}\sum_k (F_{A_k}(\cdot)-1)\qquad
\mbox{\rm in}\quad L^2(\Om,\mathbb P),
\ee
where $A_1,\ldots,A_k$ is a measurable partition of $X$, and each $F_{A_k}$
is assumed to be normalized so that $\E F_{A_k}(\cdot)=1$.
\end{theorem}

In Sect.~\ref{sect:log}, we apply this theorem to compute the spaces of additive
functionals for the Gaussian and Poisson processes. 

Note that multiplicative and additive functionals can be defined in an obvious
way for general processes. However, if the process is not a process with independent 
values, then it does not determine a factorization, so that formula~(\ref{log})
does not make sense, and the logarithm is not defined.

Thus, if we have two factorized measure spaces with isomorphic Hilbert
factorizations of the corresponding $L^2$ spaces,
then this isomorphism can be determined by indicating
the correspondence either between all chaoses, either
between additive functionals, or between multiplicative
functionals. We will construct all
these correspondences for the Poisson--Gauss and L\'evy--Gauss
isomorphisms.

\subsection{General L\'evy processes}
\label{sect:levy}

The notion of a generalized random process was developed by I.~M.~Gelfand
\cite{Gel55} (see also \cite{GelVil61})
and K~It\^o \cite{Ito54}. Generalized processes with independent values 
are considered in detail in
\cite{GelVil61}, see also \cite{Feld71}. A standard reference on 
L\'evy processes in $\mathbb R^n$ is the monograph
\cite{Bert96}. 

Let $(X, \nu)$ be a standard Borel space with a continuous finite measure
$\nu$. Denote by $\cal F$ the linear space of ($\bmod 0$ classes of)
bounded measurable functions
on $X$. A generalized random process\footnote{Recall that a generalized
random process on a real topological vector locally convex space
$\cal L$ is a continuous linear map
$a\mapsto\langle a, \cdot\rangle$
from $\cal L$ to the space $L_0(\Omega, \frak A, \mathbb P)$
of random variables on a probability space
$(\Omega, \frak A, \mathbb P)$. A generalized random process induces 
a weak distribution (a measure on cylinder sets) on the conjugate space
$L^*$, which sometimes can be extended to a probability measure.} 
on the space $\cal F$ is called a
{\it process with independent values (L\'evy process)} if for any
functions $a_1,a_2\in \cal F$ such that $a_1(x)a_2(x)=0$ a.e.,
the random variables $\langle a_1,\cdot\rangle$ and $\langle a_1,\cdot\rangle$
are independent. A process is called {\it homogeneous} if it is invariant
under measure-preserving transformations of the space $(X, \nu)$.

It is well-known (see, e.g., \cite{Feld71}) that homogeneous processes
with independent values are described by the L\'evy--Khintchin theorem.
For each homogeneous process with independent values on
$(X, \nu)$, there exists a Borel measure
$\Pi$ on $\mathbb R$ such that $\Pi(\{0\})=0$ and
$\int_\mathbb R\frac{t^2}{1+t^2}d\Pi(t)<\infty$ (the L\'evy--Khintchin measure),
a nonnegative number $\si^2\in\mathbb R_+$ (the Gaussian variance), and
a number $c\in\mathbb R$ (the drift) such that for every function
$a\in\cal F$,
$$
\E e^{i\langle a,\cdot\rangle}=
\exp\left(\int_X\log\phi(a(x))d\nu(x)\right),
$$
where
$$
\log\phi(y)=ic y-\frac{\si^2y^2}{2}+
\int_{\mathbb R}\left(e^{ity}-1-\frac{ity}{1+t^2}\right)d\Pi(t),
$$
and the parameters $\Pi$, $\si^2$, and $c$ are uniquely determined by the 
process. In particular, $\Pi$ is the measure of jumps of the process.

The following lemma is well-known.

\begin{lemma}\label{l:levyfact}
Let $\eta$ be an arbitrary process with independent values 
on a standard Borel space $(X,\frak B)$ with a continuous finite measure $\nu$. 
For each $A\in\frak B$,
let $\zeta_\eta(A)$ be the $\si$-field generated by the restriction
of $\eta$ to $A$. Then $\zeta_\eta$ is
a measure factorization over the Borel space $(X,\frak B)$,
which can be extended to a measure
factorization over the algebra of $\bmod0$ classes of measurable functions
(i.e., over the corresponding Lebesgue space).
\end{lemma}

According to the general construction of Lemma~\ref{l:meas2hilb}, 
a process $\eta$ with independent values defines also
a Hilbert factorization $\xi_\eta$
in the space of square integrable functionals ${\cal H}=L^2(\eta)$
of the process; namely, 
$\xi_\eta(A)={\cal B}(H_A)\otimes 1_{H_{X\setminus A}}$, where
$H_A$ is the set of square integrable functionals that depend
only on the restriction of $\eta$ to $A$.

A very important problem is to classify processes according to
the factorizations they generate.
It turns out that all processes with 
independent values generate Fock factorizations, and the complete invariant
of this Fock factorization is the number of points in 
the support $\supp\Pi$ of the L\'evy--Khintchin measure $\Pi$
(Theorem~\ref{th:indval}).

Note that the measure factorizations generated by processes with independent values 
with the same cardinality of the L\'evy--Khintchin measure need not be isomorphic.
Feldman \cite{Feld71} showed
that {\it the measure factorizations generated by 
processes with independent values without Gaussian component
are isomorphic if and only if the measure spaces $(\mathbb R,\Pi_1)$ 
and $(\mathbb R,\Pi_2)$,
where $\Pi_1$ and $\Pi_2$ are the corresponding L\'evy--Khintchin measures, 
are isomorphic}. The Gaussian factorization is not isomorphic to
a non-Gaussian one.

The Gaussian and Poisson processes described above are processes with independent
values ($\Pi=0$, $c=0$  and $\si=0$, $\Pi=\de_1$,
$c=\frac12$, respectively\footnote{It is sometimes convenient to think that in the 
Gaussian case the L\'evy--Khintchin measure 
$\Pi$ is concentrated at 0.}). Moreover, each process
$\eta$ with independent values can be uniquely decomposed into the sum 
$\eta=\eta_0+\eta_1+\eta_2$, where $\eta_0$ is a deterministic component: 
$\langle\eta_0,a\rangle=c\cdot\int_Xa(x)d\nu(x)$, and
$\eta_1$ and $\eta_2$ are {\it independent} processes with independent
values, $\eta_1$ being a Gaussian process, and
$\eta_2$ being a purely jump (i.e., having no
Gaussian component) process
with zero drift; thus, in order
to study the spaces of functionals of a L\'evy process,
it suffices to consider separately Gaussian processes and purely jump 
processes without drift. The Gaussian case was considered above.
If $\eta$ is a L\'evy process without Gaussian component and drift, then
it can be uniquely recovered from the measure of jumps, which is the Poisson
process $\pi_{\nu\times\Pi}$ on the space $X\times\mathbb R$ with the product
mean measure $\nu\times\Pi$. Thus the space
$L^2(\eta)$ of square integrable functionals of a L\'evy process
without Gaussian component can be identified with
the space $L^2(\pi_{\nu\times\Pi})$ of square integrable functionals of 
the Poisson process on the space $(X\times\mathbb R, \nu\times\Pi)$.
Clearly, this identification preserves the structure of a factorization
over $(X,\nu)$. Thus the
{\em study of the factorizations generated by an arbitrary process
with independent values reduces to the case of the Gaussian and
Poisson processes}.

\smallskip\noindent{\bf Remark.} In the case when the 
Gaussian component and drift are zero, and the L\'evy--Khintchin measure is concentrated
on $\mathbb R_+$ and satisfies the condition
\be
\label{subord}
\int_0^\infty(1-e^{-s})d\Pi(s)<\infty,
\ee
it is more convenient to define the L\'evy process by the Laplace transform
$$
\E e^{-\langle a,\eta\rangle}=
\exp\left(\int_X\log\psi_\Pi(a(x))d\nu(x)\right),
\qquad a\ge0,
$$
where $\psi_\Pi$ is the Laplace transform of the infinitely divisible
distribution $F_\Pi$ with the L\'evy--Khintchin measure $\Pi$:
$$
\psi_\Pi(t)=\exp\left(-\int_0^\infty(1-e^{-ts})d\Pi(s)\right).
$$
If $X=\mathbb R$, then such processes are
{\it subordinators}\footnote{Recall (see, e.g., \cite{Bert96, Ki93}) that
a subordinator is a homogeneous process on $\mathbb R$
with independent positive increments.}, thus a L\'evy
processes that satisfies the above conditions will be called
a {\it generalized subordinator}.

In the case of subordinators, the explicit construction of the process by means of
the corresponding  Poisson process looks as follows
(see, e.g., \cite[Ch.~8]{Ki93}).
Consider a Poisson point process on the space
$X\times\Bbb R_+$ with the mean measure $\nu\times\Pi$.
We associate with a realization $\pi=\{(x_i,z_i)\}$ of this process the 
measure
\begin{equation}
\eta=\sum_{(x_i,z_i)\in\pi} z_i\de_{x_i}.
\label{poisson}
\end{equation}
Then $\eta$ is the generalized subordinator with the L\'evy measure $\La$.
In particular,
it follows that the distribution of a generalized subordinator is concentrated on the cone
$D^+$, where
$$
D=\left\{\sum z_i\de_{x_i},\;x_i\in X,\,z_i\in\Bbb R,\sum|z_i|<\infty \right\}
$$
is the real linear space of all finite real discrete measures
on $X$, and 
$D^+=\{\sum z_i\de_{x_i}\in D:\;z_i>0\}\subset D$ is the cone in $D$ consisting
of all positive measures.

Note that in the case of subordinators, it is easy to present
explicitly the isomorphism of factorizations from Feldman's theorem (see 
above in this section). Namely, consider two generalized subordinators
$\eta_1$ and $\eta_2$ with L\'evy measures $\Pi_1$ and $\Pi_2$,
respectively, and let
$T:\mathbb R_+\to \mathbb R_+$ be the map that realizes the isomorphism
of measure spaces
$(\mathbb R,\Pi_1)$ and $(\mathbb R,\Pi_2)$. Then the map
$$
T(\sum_{(x_i,z_i)\in\pi} z_i\de_{x_i})=\sum_{(x_i,z_i)\in\pi} T(z_i)\de_{x_i}
$$
realizes the isomorphism of the measure factorizations generated
by the generalized subordinators
$\eta_1$ and $\eta_2$.

\subsection{On groups of 
automorphisms of factorizations}

According to Definition~\ref{def:isomfact},
the group $\AUT(\cal H,\xi)$
of automorphisms of a factorized Hilbert space $({\cal H},\xi)$
over a Boolean algebra ${\cal A}$ consists of isometries $T$ of the space $\cal H$
such that $\bar T\circ\xi=\xi\circ S$ for some automorphism $S$ of 
the Boolean algebra $\cal A$
(recall that $\bar T$ stands for the operator in the lattice ${\cal R}(\cal H)$
generated by $T$).
The subgroup $\SAUT(\cal H,\xi)$ of {\em special} automorphisms 
consists of automorphisms that leave the base unchanged:
$\bar T\xi(A)=\xi(A)$ for all $A\in\cal A$ (i.e., $S$ is the identity 
automorphism of $\cal A$).
The subgroup $\SAUT_\psi({\cal H},\xi)\subset\SAUT(\cal H,\xi)$
consists of special automorphisms of $({\cal H},\xi)$ that leave a 
multiplicative vector $\psi\in\cal H$ (the vacuum) unchanged.
Thus there is the following natural hierarchy of groups of 
automorphisms of a factorized Hilbert space
$({\cal H},\xi)$:
$$
\SAUT_\psi({\cal H},\xi)\subset\SAUT(\cal H,\xi)\subset\AUT(\cal H,\xi).
$$
Denote by $\Aut(X,\nu)$ the group of automorphisms of the measure space
$(X,\nu)$.

For a Fock factorization, the groups $\SAUT(\cal H,\xi)$ and
$\SAUT_\psi({\cal H},\xi)$ were computed in \cite{Ar70}.

\begin{proposition}\label{prop:aut}
Let $\xi$ be a homogeneous Fock factorization in the space
${\cal H}=\EXP L^2((X,\nu);H)$.
The group $\AUT({\cal H})$ 
consists of operators of the form
\be\label{AUT}
\EXP h\mapsto e^{ib-\frac{\|\psi\|}{2}-(\psi,Uh(S^{-1}\cdot))}\cdot
\EXP(Uh(S^{-1}\cdot)+\psi),
\ee
where $b\in\mathbb R$, $\psi\in L^2((X,\nu);H)$, $S\in\Aut(X,\nu)$,
and  $U$ is a unitary operator in $L^2((X,\nu);H)$ that commutes with
all projections $P_A$ to subspaces of the form $L^2((A,\nu_A);H)$,
where $A$ is a measurable subset of $X$ and $\nu_A$ is the restriction of 
the measure $\nu$
to $A$.

The subgroup $\SAUT({\cal H})$ of special automorphisms
consists of operators (\ref{AUT}) with $S=\Id$ 
(the identity map in $X$), and the subgroup 
$\SAUT_1({\cal H})$ of vacuum-preserving special automorphisms
consists of operators (\ref{AUT}) with $S=\Id$ and $\psi=0$. 
\end{proposition}

In particular, in the case of a one-dimensional Fock factorization 
in the space ${\cal H}=\EXP L^2(X,\nu)$, it follows from
the Spectral Theorem that
$U$ is a multiplicator $h(\cdot)\mapsto a(\cdot)h(\cdot)$ by a measurable
function $a:X\to\mathbb C$ with $|a|\equiv1$. Thus $\AUT({\cal H})$
is isomorphic to the semidirect product
$(\Aut(X,\nu)\rightthreetimes \mathbb T^X)\rightthreetimes(L^2(X,\nu)\rightthreetimes \mathbb T)$,
where the semidirect product
of the additive group of the Hilbert space $L^2(X,\nu)$ with 
the unit circle $\mathbb T=\{e^{ib},\,b\in\mathbb R\}$ is
determined by the cocycle $c((\psi_1,b_1),(\psi_2,b_2))=b_1+b_2-\Im(\psi_2,\psi_1)$;
the group $\Aut(X,\nu)$ of automorphisms of the measure space $(X,\nu)$ 
acts on $\mathbb T^X=\{a:X\to\mathbb C:\,|a(x)|\equiv 1\}$ as $Sa(\cdot)=a(S^{-1}\cdot)$, 
and the pair 
$(S,a)$ sends $(\psi,b)$ to $(a(\cdot)\psi(S^{-1}\cdot),b)$. 

If a Hilbert factorization is generated by a measure factorization as
in Lemma~\ref{l:meas2hilb}, then we may consider also the group $\AUT(\zeta)$ of
automorphisms of the underlying measure factorization, which 
is obviously a subgroup of $\AUT(\cal H,\xi)$, since each automorphism of
the measure space induces an isometry in the corresponding
$L^2$ space.
This group may be different for different measure factorizations 
generating the same Hilbert factorization. 

Let we are given a measure factorization $\zeta$ of a Lebesgue space
$({\cal X}, \mu)$ with a finite measure over the base Lebesgue space
$(X, \nu)$, i.e., over the Boolean algebra of
$\bmod0$ classes of measurable sets in $(X, \nu)$. Assume that each
automorphism $T$ of the base space $(X, \nu)$ induces an automorphism
$V_T$ of the factorized space
$({\cal X}, \mu)$, in other words, the factorization is invariant
with respect to the group $\Aut(X,\nu)$. These conditions are 
satisfied if the measure factorization generates a 
{\it Fock} Hilbert factorization in $L^2({\cal X}, \mu)$ 
and, as follows from Tsirelson's theorem
(see Sect.~\ref{sect:nonfock} of the Introduction), only in this case.
Thus we have a nontrivial monomorphism
$$
\Aut(X, \nu) \to \AUT(\zeta),
$$
i.e., a monomorphic embedding of the group of all automorphisms
of a Lebesgue spaces with a finite or $\si$-finite measure into the
group of automorphisms
of a Lebesgue spaces with a finite measure; V.~A.~Rokhlin called this 
a ``dynamical system over
a dynamical system''.

Thus we obtain the problem of finding the
(operator and metric) classification
of such systems over systems that arise from the factorizations generated
by L\'evy processes. It seems that a particular case of this problem was
first considered in 1956 by K.~It\^o \cite{Ito56}. Namely, in our
terms his result can be stated as follows.
Let $(X,\nu)$ be the real line
$\mathbb R$ with the Lebesgue measure, and
let $({\cal X}, \mu)$ be the space of realizations of a L\'evy process
(understood as a process with independent values) on
$\mathbb R$. Consider the action of the one-parameter group of shifts on $\mathbb R$
on the space 
$({\cal X}, \mu)$; then for each nondegenerate L\'evy process, these actions 
are spectral isomorphic; more precisely, the corresponding one-parameter
groups always have the Lebesgue spectrum of infinite multiplicity. 
The proof in \cite{Ito56} uses arguments similar to the Wiener--It\^o
decomposition for L\'evy processes, which makes it close to our
considerations. In the same paper, a problem of metric isomorphism of
these actions for different L\'evy processes was posed. Now, using
the achievements of the ergodic theory, one can answer this question in the
affirmative.
  
\begin{theorem}
The action of the one-parameter group of shifts on the space of realizations
of any nondegenerate L\'evy process is a Bernoulli action of the group 
$\mathbb R^1$ with infinite entropy. Hence all these actions are
metrically isomorphic.
\end{theorem}
     
Note that this isomorphism does not preserve the type
of the factorization.
     
Consider the action of the whole group
$\Aut(X, \nu)$ on the space of realizations
of a L\'evy process. It is not difficult to show that the spectral
or metric isomorphism of two such actions 
implies the corresponding (Hilbert
or metric) isomorphism of the factorizations. It will follow from our
main result that the actions are spectral isomorphic for L\'evy
processes of the same dimension (i.e., with the same cardinality
of the support of the L\'evy measure) and nonisomorphic for L\'evy
processes of different dimensions. The metric
classification of actions reduces to the metric classification of
factorizations, see Sect.~\ref{sect:levy} above.

Note that the metric theory of actions of groups of automorphisms of
the Gaussian process is well developed starting from the works by
A.~N.~Kolmogorov, S.~V.~Fomin, I.~V.~Girsanov, Ya.~G.~Sinai, 
G.~Maruyama, A.~M.~Vershik, and others
(see, e.g., \cite{KSF, DynSyst}), however, for other processes with independent
values, much less is known. For example, for the gamma process, this theory
must be of interest (see, e.g., \cite{TVY01}).

\subsection{On orthogonal decompositions}
\label{sect:orthdec}

By an orthogonal decomposition in the space of square integrable
functionals of a process $\phi$ with independent values
we mean the result of the standard orthogonalization process 
in $L^2(\phi)$ applied to the symmetric tensor powers
of the subspace of additive functionals (the first chaos).
The problem of constructing such a decomposition for an arbitrary process
with independent values can be addressed within the following general scheme.

Let $\phi$ be an arbitrary process with independent values on the space $X$,
which can be regarded as a random measure on $X$.
We will construct an orthogonal decomposition in the space $L^2(\phi)$
of square integrable functionals of this process.
Set ${\cal H}_0=\mathbb C$. Let
${\cal H}_1$ be the subspace of centralized (i.e., orthogonal to constants)
{\it additive} functionals in the factorization generated by $\eta$
(if we are given the set $\cal M$ of 
square integrable multiplicative functionals in this factorization, 
then the space ${\cal H}_1$ 
can be obtained by the logarithm construction described in Theorem~\ref{th:log}).
Assuming that we have already constructed ${\cal H}_1,\ldots, {\cal H}_{n-1}$, the next space
${\cal H}_n$ is defined as the orthogonal complement to ${\cal H}_1\oplus{\ldots}\oplus{\cal H}_{n-1}$
in the subspace spanned by the functionals of order $n$, that is, by products
of $n$ additive functionals:
$$
{\cal H}_n=\overline{\left\{F_{1}(\phi)\cdot\ldots\cdot F_{n}(\phi),\;F_k\in {\cal H}_1\right\}}
\ominus({\cal H}_0\oplus\ldots\oplus {\cal H}_{n-1}).
$$
The space ${\cal H}_n$ is called the
{\it $n$th chaos} of the process $\phi$.
We have $L^2(\phi)=\oplus_{n=0}^\infty {\cal H}_n$.

Note that for general L\'evy processes, the space of {\it additive}
functionals ${\cal H}_1$ is wider than the space of {\it linear} functionals 
$L_1$.
If we apply the above orthogonalization process to $L_1$ instead of ${\cal H}_1$,
then the $n$th subspace $L_n$ will be the space of $n$-multiple 
stochastic integrals of $\phi$, but the sum of this subspaces will not exhaust
the whole space $L^2(\phi)$ (see \S\ref{sect:genisom}). 
However, for the Gaussian and Poisson processes,
the spaces of linear and additive functionals coincide, hence
${\cal H}_n=L_n$ for all $n$. 

The spaces $L_n$ can be described explicitly using the general
combinatorial approach to stochastic integrals suggested by
Rota and Wallstrom \cite{RotaW97}.
Denote by $\De_n$ the $n$th
diagonal measure of $\phi$, that is, 
$$
\De_n(A)=\phi^{\otimes n}\{(x,\ldots,x):\,x\in A\}
$$
for each measurable $A$.
Then it follows easily from the  results of \cite{RotaW97} (see also 
\cite{SegKail76}) that $L_n$ is spanned by
the $n$-multiple stochastic integrals of the form
\begin{multline}
\label{rota}
I^{(n)}_{f_1,\ldots,f_n}(\phi)=
\int f_1(x_1){\ldots}f_n(x_n)\,d\phi(x_1,{\ldots},x_n)\\=
\sum_{g\in\sn(x_1,{\ldots},x_n)}(-1)^{n-c(g)}
\prod_{(x_{i_1}\ldots x_{i_k})\in C(g)}
\int f_{i_1}(x){\ldots}f_{i_k}(x)d\De_k(x),
\end{multline}
where $\sn(x_1,{\ldots},x_n)$ is the symmetric group of degree $n$ realized
as the group of permutations of the set $\{x_1,\ldots,x_n\}$,
$C(g)$ is the set of cycles of a permutation $g$, and
$c(g)=\#C(g)$ is the number of cycles in $g$.

In Sect.~\ref{sect:orth0}, we will consider orthogonal decompositions 
for the Poisson and Gaussian processes, and in 
\S\ref{sect:genisom}, for general L\'evy processes.

\section{Canonical isomorphism between the
factorizations generated by the Gaussian and Poisson processes}
\label{sect:isom0}

In this section, we study the isometry
between the spaces of square integrable functionals of the Gaussian
and Poisson processes over the same base space.
Many authors observed
that these spaces have many common features. However,
the existence of a natural isometry between these spaces
was first established by Vershik, Gelfand, and Graev  \cite{GGV75} 
from considerations related to the representation theory of
groups of diffeomorphisms;
this isometry was explicitly
described by Neretin \cite{Ner97} as an isometry between the space
of square integrable functionals of 
the Poisson process and the so-called holomorphic model of the
boson Fock space. 

From now on  we fix the base space $(X,\nu)$, which is a 
continuous Lebesgue space. According to the general construction of
the factorization determined by a process with independent values,
the spaces
$L^2(\al)$ and $L^2(\pi)$  of square integrable functionals of
the Gaussian white noise and the Poisson process on
$(X,\nu)$, respectively, are equipped with Hilbert factorizations over
$(X,\nu)$.

\begin{theorem}
\label{th:main}
There exists a unique
unity-preserving special real\footnote{That is, sending the 
real subspace of real-valued functionals of one process to the similar
subspace for the other process.} 
isomorphism of the Hilbert factorizations
in the spaces $L^2(\al)$ and $L^2(\pi)$. The corresponding isometry
$$
\Phi :L^2(\al)\to L^2(\pi)
\label{canisom}
$$
of the Hilbert spaces is given by the following formula
on the set of multiplicative functionals: for each 
$h\in L^2(X,\nu)\cap L^1(X,\nu)$,
\be
\Phi:\;\; e^{<h,\cdot>-\frac{\|h\|^2}{2}}\mapsto
\prod_{x\in \om}(1+h(x))\cdot
e^{-\int h(x)d\nu(x)},\qquad \om\in\cal E.
\label{isom0}
\ee
\end{theorem}

\noindent{\bf Remarks. 1.} Since the map
$L^2(X,\nu)\ni h\mapsto\Phi h\in L^2(\pi)$, given by
(\ref{isom0}), is continuous in the topology of $L^2(X,\nu)$
(since the norm of the functional determined by the right-hand side equals,
as can be easily seen, $e^{\|h\|^2}$), it extends by continuity
to the whole space $L^2(X,\nu)$.

\smallskip\noindent{\bf 2.} Formula (\ref{isom0}) itself appeared,
e.g., in \cite{DeMe}, however, it was not apparently observed that
this formula determines an isometry of Hilbert spaces and
an isomorphism of factorizations.

\begin{proof} Formula (\ref{isom0}) follows from the 
formula for the boson--Poisson correspondence described in
\cite{Ner97} in terms of the holomorphic model of the Fock space
and formula
(\ref{fockgauss}) for the isomorphism between $L^2(\al)$ and
the boson Fock space. On the other hand,
it is not difficult to check directly that the map defined by (\ref{isom0})
is an isomorphism with desired properties. 

Obviously, any two isomorphisms with desired properties differ from each
other by a vacuum-preserving special
automorphism of the Fock factorization in $\EXP L^2(X,\nu)$,
i.e., by an element of the group $\SAUT_1(\EXP L^2(X,\nu))$.
Thus the uniqueness part follows from Proposition~\ref{prop:aut}.
\end{proof}

{\bf Remark.} In particular, we obtain that
the set of square integrable multiplicative functionals 
of the Poisson process
is $\{c\cdot\prod_{x\in\om}(1+h(x)),\;h\in L^2(X,\nu)\}$.
\smallskip

Our purpose is to study the isometry (\ref{isom0}) in more detail. In particular,
we would like to find its \emph{kernel},
that is, a (generalized) function $K(\om,f)$ on ${\cal E}\times \hat H$
such that for every $F\in L^2(\pi)$, 
\be
\label{kerneldef}
(\Phi^{-1} F)(\cdot)=\int_{\cal E} K(\om,\cdot)F(\om)d{\cal P}(\om).
\ee

Let $\eta$ be a process with independent values on $(X,\nu)$.   
A functional $F\in L^2(\eta)$ is called {\em singly generated}
if it depends only on the integral $\langle\eta,1\rangle$
of the process $\eta$ over the whole space $X$.
Similarly, $F$ is called {\em finitely generated} ($n$-generated)
if there is a finite measurable partition $X=A_1\cup\ldots\cup A_n$ 
of the space $X$ such
that $F$ depends only on the integrals $\langle\eta,\chi_{A_1}\rangle,\ldots,
\langle\eta,\chi_{A_n}\rangle$
of the process $\eta$ over the subsets $A_1,\ldots,A_n$.
Note that the space $L^2(\eta)$ is the projective limit 
of the subspaces of finitely generated functionals 
with respect to refinement of partitions.
Clearly, each isomorphism of the factorizations generated by
two processes $\eta_1$ and $\eta_2$ must send the set
of $n$-generated functionals of $\eta_1$ to the same set for
$\eta_2$. 

\subsection{The restriction of the canonical isomorphism to the subspace of 
singly generated functionals}
\label{sect:onepoint}

Consider the restriction of the isomorphism (\ref{isom0})
to the subspaces of singly generated functionals (see the definition above).
If $\eta$ is the Gaussian white noise on the space
$(X,\nu)$ with $\nu(X)=a$, then 
$\eta(X)=\langle\eta,1\rangle$ is the Gaussian random variable
with zero mean and variance $a$, so it is natural to identify
the space of singly generated functionals
of $\eta$ with the space
$L^2({\mathbb R}, N(0,a))$ of square integrable functions with
respect to the normal distribution (which can be also regarded
as the $L^2$ space over the Gaussian process on the space
$X=\{x\}$ that consists of a single point of weight
$a$). Similarly, the space of singly generated functionals of the Poisson process
$\pi$ on $(X,\nu)$ is identified with
the space $L^2(\mathbb Z_+, P_a)$
of sequences
$b=\{b_n\}_{n\ge0}$ with the scalar product
$$
(b,b^\prime)=e^{-a}\sum_{n=0}^\infty \frac{a^n}{n!}b_nb^\prime_n
$$
(i.e., the $L^2$ space over the Poisson process on a single-point
space).
Consider this situation in more detail.

Formula~(\ref{isom0}) takes the following form: for all $t\in\mathbb R$,
\be
e^{tx-\frac{at^2}{2}}\leftrightarrow
\left\{ e^{-at}(1+t)^k\right\}_{k=0}^\infty.
\label{isom1d}
\ee
Note that the left-hand side of (\ref{isom1d}) is the generating function
for the Hermite polynomials $H_n^a(x)$ and
the right-hand side is the generating function for the Charlier polynomials
$C_n^a(k)$ (see Appendix D). These
polynomials constitute orthogonal families in $L^2(\mathbb R,N(0,a))$ and
$L^2(\mathbb Z_+,P_a)$, respectively, and 
$\|H^a_n\|_{L^2(\mathbb R,N(0,a))}^2=\|C^a_n\|_{L^2(\mathbb Z_+,P_a)}^2=a^nn!$, 
hence
\be
\label{poly1d}
\Phi(H^a_n)=C^a_n.
\ee
We see that formula~(\ref{isom1d}) is precisely the expression of (\ref{poly1d})
in terms of generating functions.

\begin{proposition}
The kernel {\rm(\ref{kerneldef})} of the unitary isomorphism between
$L^2(\mathbb R,N(0,a))$ and $L^2(\mathbb Z_+,P_a)$ is given by the
formula
\be
\label{kernel1d}
K^a(k,x)=e^{-\frac a2-x}\frac{H_k^a(x+2a)}{a^k}.
\ee
\end{proposition}
\begin{proof}
Observe that if $\nu_1$ and $\nu_2$ are arbitrary measures on $\mathbb R$ with finite moments,
and $P^{(j)}_n$ are the orthonormalized polynomials with respect
to $\nu_j$, $j=1,2$, then
the kernel of the unitary isomorphism between
$L^2(\mathbb R,\nu_1)$ and $L^2(\mathbb R,\nu_2)$ is given by the formula
$$
K(x,y)=\sum_{n=0}^\infty P^{(1)}_n(x)P^{(2)}_n(y),
$$
provided that the series converges in $L^2(\mathbb R,\nu_1)$
for a.e.~$y$.

Thus let us consider the series 
$$
\sum_{n=0}^\infty \frac{H^a_n(x)C^a_n(k)}{a^nn!}.
$$
We would like to prove that this series converges in $L^2(\mathbb R,N(0,a))$
for each $k\in\mathbb Z_+$.
Since $H^a_n$ are orthogonal polynomials in $L^2(\mathbb R,N(0,a))$, and
$\|H^a_n\|^2={n!}{a^n}$, 
it suffices to check that 
\be
\label{1}
\sum_{n=0}^\infty\frac{|C^a_n(k)|^2}{a^nn!}<\infty.
\ee
But in view of (\ref{chform}), 
$\frac{|C^a_n(k)|^2}{a^n}=|C^a_k(n)|^2a^{n-2k}\le \text{const}\cdot n^{2k}a^{n-2k}$, 
since
$C_k^a$ is a polynomial of degree $k$, and (\ref{1}) follows immediately.

Thus we have
$$
K(k,x):=K^a(k,x)=\sum_{n=0}^\infty \frac{H^a_n(x)C^a_n(k)}{a^nn!}.
$$
Consider the generating function
$$
{\cal K}(t,x)=\sum_{k=0}^\infty\frac{t^k}{k!}K(k,x).
$$
Changing the order of summation yields
$$
{\cal K}(t,x)=\sum_{n=0}^\infty \frac{H_n^a(x)}{n!}
\sum_{k=0}^\infty\frac{t^k}{k!}\frac{C_n^a(k)}{a^n}.
$$
In view of (\ref{chform}) and (\ref{chargen}), the internal sum equals
$$
\sum_{k=0}^\infty(-1)^{n+k}\frac{t^k}{a^kk!}C_k^a(n)=(-1)^ne^t\left(1-\frac ta\right)^n.
$$
Thus
$$
{\cal K}(t,x)=\sum_{n=0}^\infty\frac{H^a_n(x)a^n}{n!}\left(\frac ta-1\right)^ne^t=
\exp\left(2t+\frac{tx}{a}-x-\frac{t^2}{2a}-\frac a2\right),
$$
which coincides, in view of~(\ref{hermgen}), with the generating function for
$e^{-\frac a2-x}\frac{H_k^a(x+2a)}{a^k}$.
\end{proof}

\smallskip\noindent{\bf Remark.} One can also check
formula~(\ref{kernel1d})
directly, using a known formula (see, e.g., \cite[formula 7.377 of the fourth Russian
edition]{GradRyzh})
$$
\frac{1}{\sqrt{2\pi}}\int_{-\infty}^{\infty}e^{-\frac{x^2}{2}}
H_n(x+a)H_m(x+b)dx=m!(-b)^{n-m}L^{n-m}_{m}(-ab),\quad m\le n,
$$
where $L^\al_n(x)$ is the Laguerre polynomial with parameter $\al$,
and the formula $C_n^a(x)=n!L^{x-n}_n(a)$ relating the Charlier and
Laguerre orthogonal polynomials.

\begin{corollary} We have proved the following identity
relating the Hermite and Charlier orthogonal polynomials
\be
\label{hermcharid}
\sum_{n=0}^\infty\frac{H^a_n(x)C^a_n(k)a^n}{n!}=
e^{-\frac a2-x}\frac{H_k^a(x+2a)}{a^k}.
\ee
In particular, denoting $C_n(\cdot)=C^1_n(\cdot)$, we obtain
\be
\sum_{n=0}^\infty\frac{H_n(x)C_n(k)}{n!}=e^{-\frac 12-x}H_k(x+2).
\ee
\end{corollary}

Thus in the case of a single-point space $X$, we have the description of the 
Gauss--Poisson isomorphism at three levels: 
\begin{itemize}
\item correspondence of multiplicative functionals: (\ref{isom1d});
\item correspondence of orthogonal polynomials: (\ref{poly1d});
\item explicit kernel: (\ref{kernel1d}).
\end{itemize}
Our purpose is to obtain the formulas of the second and third level
for general spaces.

\subsection{The restriction of the canonical isomorphism to the subspace of
finitely generated functionals}
\label{sect:finpoint}
Let  $X=A_1\cup\ldots\cup A_m$ be a measurable partition 
of the space $X$ with $\nu(A_j)=a_j$. Consider the corresponding
subspaces of finitely generated functionals in
$L^2(\al)$ and $L^2(\pi)$. Like in the case of singly generated functionals,
these subspaces can be identified with the $L^2$ spaces over the
corresponding processes on a
finite space $X=\{s_1,{\ldots},s_m\}$ that consists of $m$ points with weights
$\nu(s_j)=a_j$, $j=1,{\ldots} ,m$. Using the results of the
previous section, we obtain that
\begin{eqnarray*}
L^2(\al)&=&\bigotimes_{j=1}^m L^2(\mathbb R, N(0,a_j)),\\
L^2(\pi)&=&\bigotimes_{j=1}^m L^2(\mathbb Z_+, P_{a_j}),\\
\Phi\left(\prod_{j=1}^mH^{a_j}_{k_j}(\cdot)\right)&=&
\prod_{j=1}^m C^{a_j}_{k_j}(\cdot),
\end{eqnarray*}
and for $k=(k_1,{\ldots} ,k_m)\in \mathbb Z_+^n$ and
$x=(x_1,{\ldots} ,x_m)\in \mathbb R^n$,
\be
K(k,x)=\prod_{j=1}^me^{-\frac{a_j}{2}-x_j}\frac{H^{a_j}_{k_j}(x_j+2a_j)}{a_j^{k_j}}.
\label{kernelfindim}
\ee

\subsection{Logarithm}
\label{sect:log}

In this section, we apply the logarithmic construction described in 
Theorem~\ref{th:log}
to computing the sets of additive functionals for the Gaussian and Poisson
processes. Note that the canonical isomorphism
$\Phi$ sends the logarithmic operation in the Poissonian factorization
to the logarithmic operation in the Gaussian factorization,
and hence sends additive functionals to additive functionals.

\subsubsection*{Poisson process}
In this case
$(\Om,\mathbb P)=({\cal E}(X),\cal P)$,
and normalized multiplicative functionals are given by the formula
$$
F_A(\om)=\prod_{x\in\om\cap A}(1+h(x))e^{-\int_Ah(x)d\nu(x)},
\qquad h\in L^2(X,\nu),
$$
for each measurable subset $A\subset X$.

\begin{lemma}
\label{l:logpois}
\quad $\LOG F_A(\om)=\sum_{x\in\om\cap A}h(x)-\int_Ah(x)d\nu(x)$.
\end{lemma}
\begin{proof}
Denote the right-hand side by $G_A(\om)$.
Note that
\begin{eqnarray*}
\|\sum_k (F_{A_k}(\cdot)-1)-G(\cdot)\|^2&=&
\|\sum_k (F_{A_k}(\cdot)-1-G_{A_k}(\cdot))\|^2\\&=&
\sum_k\|F_{A_k}(\cdot)-1-G_{A_k}(\cdot)\|^2,
\end{eqnarray*}
since the restrictions of the Poisson
process to disjoint subsets are independent.
Using Campbell's theorem for sums and products over Poisson processes
(see, e.g., \cite{Ki93}, Sect.~3.2, 3.3),
it is not difficult to compute that in our case
$$
m_F(A)=\int_{A}h^2(x)d\nu(x), 
$$
and
$$
\|F_{A_k}(\cdot)-1-G_{A_k}(\cdot)\|^2=
e^{m_F(A_k)}-1-m_F(A_k).
$$
Assume that
$m_F(A_k)<\de$ for all $k$. Then
for sufficiently small $\de$ we have
$$
\|F_{A_k}(\cdot)-1-G_{A_k}(\cdot)\|^2< C\cdot m_F(A_k)^2.
$$
Hence
\begin{eqnarray*}
\|\sum_k (F_{A_k}(\cdot)-1-G_{A_k}(\cdot))\|^2\le 
C\cdot\sum_k m_F(A_k)^2\\
\le C\cdot\sum_k\de\cdot\int_{A_k}h^2(x)d\nu(x)\le C\cdot\|h\|^2\cdot\de,
\end{eqnarray*}
which is arbitrarily small for sufficiently small $\de$, and we are done.
\end{proof}

Thus we obtain that the space of additive functionals (``the first chaos'')
of the Poisson process is 
\be
\label{poislin}
\left\{\sum_{x\in\om}h(x)-\int_Xh(x)d\nu(x), \quad h\in L^2(X,\nu)\right\}.
\ee

\subsubsection*{Gaussian process}
Let us now compute the logarithm for the Gaussian processes.
In this case a normalized multiplicative functional is given by the formula
$$
F_A(\eta)=e^{\langle h,\eta_A\rangle-\frac12\int_Ah^2(x)d\nu(x)},
\qquad h\in L^2(X,\nu),
$$
where $\eta_A$ is the restriction of $\eta$ to $A\subset X$.

\begin{lemma}
\label{l:loggauss}
\quad
$\LOG F(\eta)=\langle h,\eta\rangle$.
\end{lemma}
\begin{proof}
Let $G_A(\eta)=\langle h,\eta_A\rangle$.
It is easy to compute that in this case
$$
m_F(A)=\log\E\|e^{\langle h,\eta_A\rangle-\frac12\int_Ah^2(x)d\nu(x)}\|^2
=\int_Ah^2(x)d\nu(x),
$$
and
\begin{eqnarray*}
\|F_{A}(\cdot)-1-G_{A}(\cdot)\|^2&=&
\E\|e^{\langle h,\eta_A\rangle-\frac12\int_Ah^2(x)d\nu(x)}-1-
\langle h,\eta_A\rangle\|^2\\
&=&e^{m_F(A)}-1-m_F(A),
\end{eqnarray*}
exactly as in the Poissonian case,
so the proof just reproduces the proof of Lemma~\ref{l:logpois}. 
\end{proof}

\subsection{Correspondence of orthogonal decompositions (chaoses)}
\label{sect:orth0}

In this section, we will state the canonical isomorphism
(\ref{isom0}) in terms of orthogonal decompositions.
Recall that by the orthogonal decomposition we mean the result
of the standard orthogonalization process in the Hilbert space
$L^2$ applied to the symmetric tensor powers of the 
subspace of additive functionals (the first  chaos). The general scheme
for constructing such decomposition is described in
Sect.~\ref{sect:orthdec}. Recall also that for the Gaussian and Poisson
processes, the spaces of additive and linear functionals coincide,
hence the space of the $n$th chaos ${\cal H}_n$ coincides with the space
of $n$-multiple stochastic integrals $L_n$, and the construction of
the orthogonal decomposition can be performed using the combinatorial scheme
described in Sect.~\ref{sect:orthdec}.

The orthogonal decomposition
\be
L^2(\al)=\bigoplus_{n=0}^\infty {\cal H}_n
\label{gaussdec}
\ee
in the space of square integrable functionals of the Gaussian white noise
is the well-known  Wiener--It\^o--Cameron--Martin decomposition.
Though the corresponding formulas are classical, it is instructive to 
observe how they can be obtained in the general combinatorial scheme.
If $\phi$ is the Gaussian white noise on $(X,\nu)$,
then $\De_2=\nu$, $\De_3=\De_4=\ldots=0$ (\cite{RotaW97}, Example G), whence
${\cal H}_n=L_n$ is spanned by the functionals of the form
\be
\label{multherm}
{\frak H}^{(n)}_{f_1,{\ldots},f_n}(\cdot)=\sum_{g\in\inv_n}
\prod_{i\in C_1(g)}\langle f_i,\cdot\rangle
\prod_{\{j,k\}\in C_2(g)} \left(-(f_j,f_k)\right),
\ee
where
$\inv_n$ is the set of all involutions in $\sn$, $C_k(g)$ is the number of cycles
of length $k$ of a permutation $g\in\sn$, and
$(f_j,f_k)=\int_X f_j(x)f_k(x)d\nu(x)$ is the scalar
product in $L^2(X,\nu)$. 

\begin{definition}
\label{def:hermite}
The functional ${\frak H}^{(n)}_{f_1,{\ldots},f_n}$  is called
the $n$th {\it generalized Hermite functional}. 
\end{definition}

In particular, for $f_1=f_2={\ldots} =f_n=f$, we obtain (see Appendix D,
(\ref{hermcomb}))
\be
{\frak H}^{(n)}_{f,{\ldots},f}(\cdot)=\sum_{g\in\inv_n}
\langle f,\cdot\rangle^{c_1(g)}  \cdot
(-\|f\|^2)^{c_2(g)}=H_n^{\si}(\langle f,\cdot\rangle),
\ee
that is, $I_n^{g}(f,{\ldots},f)$ is
the $n$th ordinary Hermite polynomial in $\langle f,\cdot\rangle$ with parameter $\si=\|f\|^2$.

Note that in terms of the Fock space $\EXP H$, the subspace ${\cal H}_n$ is 
precisely the $n$-particle subspace $S^nH$.

The corresponding orthogonal decomposition for the Poisson process
was first discussed by It\^o \cite{Ito51}
and explicitly constructed by Ogura \cite{Og72}. Within the combinatorial approach of
\cite{RotaW97}, it is obtained as follows. In this case $\phi$ is the
centralized  Poisson process on $(X,\nu)$, that is,
$\phi=\sum_{x\in\om}\de_x-\nu$, where $\om$ is the Poisson process on $(X,\nu)$,
and the diagonal measures equal $\De_2=\De_3=\ldots=\om$ 
(\cite{RotaW97}, Example CP). 
Thus we have
\be
\label{poisdec}
L^2({\cal E},{\cal P})=\bigoplus_{n=0}^\infty V_n,
\ee
where $V_0=\mathbb C$ and $V_n$ is spanned by the $n$-multiple
stochastic integrals given by the formula
\begin{multline}
\label{multchar}
\frak C^{(n)}_{f_1,{\ldots},f_n}=\sum_{g\in\sn}(-1)^{n-c(g)}
\prod_{i\in C_1(g)}\left(\sum_{x\in\om}f_i(x)-\int_Xf_i(x)d\nu(x)\right) \\
\cdot
\prod_{(x_{i_1},{\ldots},x_{i_k})\in C(g)}
\sum_{x\in\om}f_{i_1}(x){\ldots}f_{i_k}(x).
\end{multline}

\begin{definition}
\label{def:charlier}
The functional $\frak C^{(n)}_{f_1,{\ldots},f_n}$  is called
the $n$th {\it generalized Charlier functional}. 
\end{definition}

If $f_1={\ldots} =f_n=\chi_A$, where
$\chi_A$ is the characteristic function of a measurable set $A\subset X$,
then
%\begin{multline}
\be
\frak C^{(n)}_{\chi_A,{\ldots},\chi_A}
=\sum_{g\in\sn}
\left(\#(\om\cap A)-\nu(A)\right)^{c_1(g)} \cdot
\#(\om\cap A)^{-c_2(g)+c_3(g)-\ldots}=C_n^\si(\#(\om\cap A)),
%\end{multline}
\ee
that is, $\frak C^{(n)}_{\chi_A,{\ldots},\chi_A}$ is 
the $n$th ordinary Charlier polynomial with parameter $\si=\nu(A)$
(see Appendix D, (\ref{charcomb})). 
(Note that in the Poissonian case, unlike the
Gaussian one, the functional $\frak C^{(n)}_{f,{\ldots},f}$ 
with an arbitrary function $f$ is not
an ordinary Charlier polynomial. The reason is that 
all linear functionals of the Gaussian process have Gaussian distributions,
while in the Poissonian case only integrals over subsets of $X$
have Poisson distributions.)
In particular, it follows from (\ref{multchar}) that
the first chaos of the Poisson process consists of functionals of the form
$\frak C^{(1)}_{f}(\pi)=\sum_{x\in\pi}f(x)-\int_Xf(x)d\nu(x)$
(cf.~(\ref{poislin})),
and the second chaos is generated by the functionals of the form
\be
\label{poisquadr}
\frak C^{(2)}_{f,g}(\pi)=\frak C^{(1)}_{f}(\pi)\frak C^{(1)}_{g}(\pi)-
\sum_{x\in\om} f(x)g(x).
\ee

\begin{corollary}
The canonical isomorphism $\Phi$ sends the generalized Hermite functional to the
corresponding generalized Charlier functional:
\be
\Phi{\frak H}^{(n)}_{f_1,{\ldots},f_n}=\frak C^{(n)}_{f_1,{\ldots},f_n}.
\ee
\end{corollary}

\subsection{Kernel}
\label{sect:kernel}
Let $T$ be an isometry of the Hilbert spaces $L^2(A,\mu)$ and
$L^2(B,\nu)$. In some cases this isometry can be represented in
the integral form
$$
(TF)(\cdot)=\int_{A}K(x,\cdot)F(x)d\mu(x),\quad F\in L^2(A,\mu),
$$
where $K$ is a (perhaps, generalized in some sense) function of two variables
on the space $A\times B$ called the kernel of the isometry. In this section, 
we will find the kernel  (\ref{kerneldef}) of the Poisson--Gauss isometry
$\Phi$, assuming, for the sake of simplicity, that
$X=[0,1]$ and $\nu$ is the Lebesgue measure on $[0,1]$, 
i.e., the kernel of the isometry between the spaces of square integrable
functionals of the standard white noise on the interval
$[0,1]$ and the homogeneous Poisson process on
$[0,1]$ with unit rate. The
case of an arbitrary continuous Lebesgue space $(X,\nu)$ 
is completely analogous. In our case the kernel turns out to be
``almost'' ordinary function, namely,
for any measurable sets $A\subset\cal E$
and $B\subset\hat H$ (recall that $\cal E$ is the set of configurations in
the space $X$, i.e., the space of realizations of the Poisson process, and
$\hat H$ is the nuclear extension of the space
$L^2(X,\nu)$, which is the space of realizations of the Gaussian white noise),
set $\rho(A,B)=\int_A\int_BK(\om,\eta)d{\cal P}(\om)d\mu(\eta)$.
Then $\rho$ is an additive set function on ${\cal E}\times\hat H$. 
Thus $K$ can be regarded as the density of a signed measure
(of infinite variation) on ${\cal E}\times\hat H$.

Let us introduce the following notation.
Given a point configuration $\om\in{\cal E}$ (since the parameter measure
$\nu$ is continuous,
all configurations of the Poisson process are simple, i.e., each point
has multiplicity one; thus we may consider only simple
(multiplicity-free) configurations),
denote by $\Pi_{\le2}(\om)$ the set
of partitions of the set $\om$ into subsets consisting of at most two points.

For example, if $\om=\{x,y,z\}$, then
$$
\Pi_{\le2}(\om)=\left\{
\{ \{x\},   \{y\}, \{z\} \},\;
\{ \{x,y\}, \{z\}        \},\;
\{ \{x\},   \{y,z\}      \},\;
\{ \{x,z\}, \{y\}        \}
\right\}.
$$

For each partition $R\in\Pi_{\le2}$, let $C_k(R)$ be the set of $k$-point
subsets in $R$, $k=1,2$, and $|R|=\#C_2(R)$.

\begin{theorem}
The kernel (\ref{kerneldef}) 
of the isomorphism (\ref{isom0}) between $L^2(\al)$ and
$L^2(\pi)$ is given by the following formula:
\be
\label{kernel}
K(\om,\eta)=e^{-\frac12-\langle\eta,1\rangle}
\sum_{R\in\Pi_{\le2}(\om)}(-1)^{|R|}\prod_{z\in C_1(R)}(\eta+2)(z)
\prod_{\{x,y\}\in C_2(R)}\de(x-y)
\ee
for almost all $\om$, $\eta$.
\label{th:kernel}
\end{theorem}
(Here $\eta+2$ is the generalized function
$\langle\eta+2,h\rangle=\langle\eta,h\rangle+2\int h(t)dt$, and
the product is the direct product of generalized functions.)

\begin{proof}
It suffices to check that (\ref{kerneldef}) holds for multiplicative functionals $F$, that is, to
show that
\be
\label{toprove}
\E\left(\prod_{x\in\om}(1+h(x))K(\om,\eta)\right)=
\exp\left(-\frac{\|h\|^2}{2}+\langle \eta+1,h\rangle\right)
\ee
for all $h\in L^2(X,\nu)$. It is more convenient to 
rewrite the sum in~(\ref{kernel}) over \emph{permutations} rather 
than partitions. Given an $n$-point configuration $\om=\{x_1,{\ldots},x_n\}$,
let $\sn(\om)$ be the symmetric group of degree $n$
realized as the group of all
permutations of the set $\{x_1,\ldots,x_n\}$. 
Let $C_i(g)$ be the set
of all cycles of length $i$ in a permutation $g\in\sn$ and set  
$c_i(g)=\#C_i(g)$. Finally, denote by $\inv(\om)$ the subset
of $\sn(\om)$ consisting of all involutions (i.e., permutations with
cycles of length at most two).
Recall that the number of points of the homogeneous Poisson process on $[0,1]$
obeys the Poisson distribution
with parameter one, and the conditional distribution of these points,
given that the number of points is equal to $n$, coincides with the distribution of
$n$ i.i.d.~variables with the uniform distribution on $[0,1]$.
Then the left-hand side of~(\ref{toprove}) equals
\begin{multline*}
e^{-3/2-\langle\eta,1\rangle}
\sum_{n=0}^\infty\frac1{n!}
\int_0^1\ldots\int_0^1(1+h(x_1))\ldots(1+h(x_n))\\ \cdot
\sum_{g\in\inv(x_1,\ldots,x_n)}
\prod_{\{x_i,x_j\}\in C_2(g)}(-\de(x_i-x_j))\cdot\prod_{x_k\in C_1(g)}(\eta(x_k)+2)\;
dx_1\ldots dx_n.
\end{multline*}
The contribution of each pair $\{x,y\}\in C_2(g)$ is equal to
$$
-\int_0^1\int_0^1(1+h(x))(1+h(y))\de(x-y)dxdy=-\int_0^1(1+h(x))^2dx,
$$
and the contribution of each element $x\in C_1(g)$ is equal to
$$
\int_0^1(1+h(x))(\eta(x)+2)dx=\langle\eta+2,1+h\rangle.
$$
Thus the sum under consideration equals
\be
\label{s1}
e^{-3/2-\langle\eta,1\rangle}
\sum_{n=0}^\infty\frac1{n!}\sum_{g\in\inv_n}t_1^{c_1(g)}t_2^{c_2(g)},
\ee
where
\begin{eqnarray}
t_1&=&\langle\eta+2,1+h\rangle,\nonumber\\
t_2&=&-\int_0^1(1+h(x))^2dx.\nonumber
\label{t}
\end{eqnarray}
But the sum in~(\ref{s1}) is just the augmented cycle index
$\tilde Z(\sn)[t_1,t_2,0,0,\ldots]$ (see Appendix C, (\ref{cyndex})).
Thus applying (\ref{cycleindex})
with $z=1$ we obtain that~(\ref{s1}) is equal to
$$
\exp\left(-3/2-\langle\eta,1\rangle+\langle\eta+2,1+h\rangle-
\frac12\int(1+h(x))^2dx\right),
$$
and~(\ref{toprove}) follows by trivial computations.
\end{proof}

\smallskip\noindent{\bf Remark.}
There is another proof of Theorem~\ref{th:kernel} which allows one to derive formula
(\ref{kernel}) rather than to check it. The idea of this proof is as follows. Observe that
$$
L^2(\al)=\varprojlim A_n,
$$
where $A_n$ is the subspace consisting of functionals $F(\eta)$ depending only
on $\langle\eta,\chi_{[0,\frac1n]}\rangle,\ldots,
\langle\eta,\chi_{[\frac{n-1}{n},1]}\rangle$.
Obviously, $A_n$ is isometric to
$$
\bigotimes_{j=1}^m L^2(\mathbb R, N(0,1/n)).
$$
Then one should apply (\ref{kernelfindim}) and pass
to the limit.

\smallskip\noindent{\bf Example 1.}
Let $A_n$ be the subset in $L^2(\pi)$ consisting of functions
supported by $n$-point configurations, $n=1,2,\ldots$. Then it follows from (\ref{kernel})
that the image of $A_1$ under the canonical isomorphism is the
subspace of functions of the form
$$
e^{-\frac32-\langle\eta,1\rangle}\langle\eta+2,f\rangle,
\quad f\in L^2([0,1],\nu);
$$
the image of $A_2$ consists
of functions of the form
$$
e^{-\frac32-\langle\eta,1\rangle}
\frac{1}{2!}\left[
\int_0^1\int_0^1f(x,y)(\eta+2)(x)(\eta+2)(y)dxdy
-\int_0^1f(x,x)dx
\right],
$$
where $f\in L^2([0,1]\times[0,1],\nu\times\nu)$; and the image of $L_3$
is
\begin{eqnarray*}
e^{-\frac32-\langle\eta,1\rangle}&\cdot&\frac{1}{3!}\left[
\int_0^1\int_0^1\int_0^1f(x,y,z)(\eta+2)(x)(\eta+2)(y)(\eta+2)(z)dxdydz\right.\\
&-&\int_0^1\int_0^1f(x,x,z)(\eta+2)(z)dxdz
-\int_0^1\int_0^1f(x,y,y)(\eta+2)(x)dxdy \\
&-&\left.\int_0^1\int_0^1f(z,y,z)(\eta+2)(y)dydz\right],
\end{eqnarray*}
$f\in L^2([0,1]\times[0,1]\times[0,1],\nu\times\nu\times\nu)$.

\smallskip\noindent
{\bf Example 2.} Each function $h\in L^2(X,\nu)$ determines a ``linear'' functional
$F_h(\om)=\sum_{x\in\om}h(x)$ of the Poisson process.
Let us compute its image in $L^2(\al)$. We have
\begin{multline*}
\E K(\om,\eta)F_h(\om)=e^{-3/2-\langle\eta,1\rangle}
\sum_{n=0}^\infty\frac1{n!}\int_0^1{\ldots} \int_0^1
\sum_{l=1}^n h(x_l)\\
\sum_{g\in\inv(x_1,\ldots,x_n)}
\prod_{\{x_i,x_j\}\in C_2(g)}(-\de(x_i-x_j))\cdot\prod_{x_k\in C_1(g)}
(\eta+2)(x_k)dx_1\ldots dx_n.
\end{multline*}
(Recall that $C_k(g)$ is the set of cycles of length 
$k$ in a permutation $g$.)
It is easy to check that each summand $h(x_l)$ contributes
$$
\begin{cases}
\langle h,\eta+2\rangle(\langle1,\eta\rangle+2)^{c_1(g)-1}(-1)^{c_2(g)},
&\text{ if }x_l\in C_1(g)\\
\langle h\rangle(\langle1,\eta\rangle+2)^{c_1(g)}(-1)^{c_2(g)},
&\text{ if }x_l\in C_2(g).
\end{cases}
$$
Thus the sum under consideration equals
$e^{-3/2-\langle\eta,1\rangle}(S_1+S_2)$, where
\begin{eqnarray*}
S_1&=&\langle h,\eta+2\rangle
\sum_{n=0}^\infty\frac1{n!}\sum_{g\in\inv(x_1,\ldots,x_n)}
c_1(g)(\langle1,\eta\rangle+2)^{c_1(g)-1}(-1)^{c_2(g)},\\
S_2&=&2\langle h\rangle
\sum_{n=0}^\infty\frac1{n!}\sum_{g\in\inv(x_1,\ldots,x_n)}
c_2(g)(\langle1,\eta\rangle+2)^{c_1(g)}(-1)^{c_2(g)}.
\end{eqnarray*}
Note that the sum in $S_1$ is the derivative of the augmented cycle index
$\ti Z(\sn)$ in $t_1$ calculated at
$t=(\langle1,\eta\rangle+2,-1,0,0,\dots)$,
hence $S_1=\langle h,\eta+2\rangle e^{\langle1,\eta\rangle+3/2}$.
Similarly, $S_2$ is the derivative of the same cycle index in $t_2$, thus
$S_2=-\langle h\rangle e^{\langle1,\eta\rangle+3/2}$,
where $\langle h\rangle=\int_Xh(x)d\nu(x)$,
and simple computations show that the image of the functional $F_h$
in $L^2(\al)$ equals
$$
F_h(\eta)=\langle h,\eta\rangle+\int_X hd\nu(x),
$$
in agreement with (\ref{poislin}).

\section{Isomorphism of the factorizations generated by
general L\'evy processes}
\label{sect:genisom}

The purpose of this section is to apply the results on the Poisson--Gauss
isomorphism to general L\'evy processes.
Recall that, as was mentioned in Sect.~\ref{sect:levy}, 
the  study of the factorizations generated by general L\'evy processes 
reduces to the study of the Poissonian and Gaussian factorizations.
We emphasize that we reduce the general case to the Poisson--Gauss
one using the universality of the isomorphism with respect to the base.

As mentioned in Sect.~\ref{sect:levy}, the space 
$L^2(\eta_\Pi)$ of square integrable functionals of a L\'evy process
$\eta_\Pi$ with L\'evy--Khintchin measure $\Pi$
can be identified with the space
$L^2(\pi_{\nu\times\Pi})$ of square integrable functionals of the Poisson process
$\pi_{\nu\times\Pi}$ on the direct product
$X\times\mathbb R$ with the mean measure $\nu\times\Pi$.
Thus it is natural to introduce the white noise
$\al_{X\times\mathbb R}$
on the space $(X\times\mathbb R,\nu\times\Pi)$.
Note that this process may be also regarded as the
$L^2(\mathbb R,\Pi)$-valued white noise $\al^{L^2(\mathbb R,\Pi)}$
on the space
$(X,\nu)$, i.e., one may identify $L^2(\al_{X\times\mathbb R})$ 
with the homogeneous Fock space
$\EXP L^2((X,\nu);L^2(\mathbb R, \Pi))$.
The spaces $L^2(\al^{L^2(\mathbb R,\Pi)})$ and $L^2(\eta_\Pi)$
are equipped with natural Hilbert factorizations over
$(X,\nu)$.

\begin{theorem}
\label{th:levy}
There exists a unity-preserving isometry (which is an
isomorphism of Hilbert factorizations)
$$
\Phi :L^2(\al^{L^2(\mathbb R,\Pi)})\to L^2(\eta_\Pi).
$$
On the set of multiplicative functionals,
it is given by the following formula: for each 
$h\in L^2(X\times\mathbb R,\nu\times\Pi)\cap  L^1(X\times\mathbb R,\nu\times\Pi)$,
\be
\Phi:\;\; e^{<h,\cdot>-\frac{\|h\|^2}{2}}\mapsto
\prod_{i}(1+h(x_i,t_i))\cdot
e^{-\int\int h(x,t)d\nu(x)d\La(t)},\quad \eta=\sum_i t_i\de_{x_i}\in D.
\label{isomlevy}
\ee
This is the unique real special vacuum-preserving automorphism of
Hilbert factorizations that acts identically on 
the space of values $L^2(\mathbb R,\Pi)$.
\end{theorem}

\noindent{\bf Remark.} When we say that the isomorphism
acts identically on the space of values, we mean that
it is an isomorphism of factorizations over the space
$(X\times\mathbb R,\nu\times\Pi)$. In other words,
consider the restriction of the isomorphism to 
the first chaos, which can be identified with
$L^2((X,\nu);L^2(\mathbb R, \Pi))$ for both processes. Then for an arbitrary
set $A\subset L^2(\mathbb R, \Pi)$, the isomorphism preserves the subset
of $L^2((X,\nu);L^2(\mathbb R, \Pi))$ that consists of functions
whose values lie in $A$.

\begin{proof}
Follows immediately from Theorem~\ref{th:main}
and the above observations.
\end{proof}

Without assuming that the isomorphism acts identically on the
space of values, the above isomorphism is not unique. Indeed,
apply Proposition \ref{prop:aut} and observe that 
the set of operators in $L^2((X,\nu);\, L^2(\mathbb R, \Pi))$
that commute with all projections $P_A$ for $A\subset X$
is $L^\infty(X,\nu)\otimes{\cal B}(L^2(\mathbb R,\Pi))$. 
Thus in this case the group $SAUT_1$  
is generated by operators of the form
$\EXP (h_1(\cdot)\otimes h_2(\cdot))\mapsto 
\EXP((a(\cdot)h_1(\cdot))\otimes(Uh_2(\cdot))),$
where $h_1\in L^2(X,\nu)$, $h_2\in L^2(\mathbb R,\Pi)$,
$a$ is a complex-valued measurable function on $X$ with $|a|\equiv 1$,
and $U$ is an arbitrary unitary operator in $L^2(\mathbb R,\Pi)$.

Theorem~\ref{th:levy} implies immediately 

\begin{theorem}
\label{th:indval}
Let $\eta$ be a L\'evy process on the space $(X,\nu)$
with L\'evy measure $\Pi$.
Then the Hilbert factorization determined by $\eta$ is a
homogeneous Fock factorization, and its dimension
is equal to the number of points in $\supp \La$.
\end{theorem}

Theorem~\ref{th:levy} states the L\'evy--Gauss isomorphism at the level of
multiplicative functionals. Let us now 
describe it in terms of orthogonal decompositions, applying 
the general scheme described in
Sect.~\ref{sect:orthdec} and assuming for simplicity that the process
is a generalized subordinator.
By (\ref{poisdec}) we have
$$
L^2(\eta_\Pi)=\bigoplus_{n=0}^\infty V_n,
$$
where $V_n$ is the space generated by the generalized Charlier
functionals of order $n$ (\ref{multchar}) 
in the space $(X\times\mathbb R,\nu\times\Pi)$.
In particular, additive functionals of
the L\`evy process are of the form
$$
\sum_i h(x_i,t_i),\qquad h\in L^2(X\times\mathbb R_+,\nu\times\La).
$$
As we have seen above, in the case of Gaussian and Poisson processes,
all additive functionals are linear. However, in the case when 
$\supp\La$ consists of more than one point, that is,
the L\'evy process is neither Gaussian nor Poisson,
{\it the space of additive functionals does not coincide
with the space of linear functionals}, which are given by 
\be
\langle a,\eta\rangle=\int_X a(x)d\eta(x)=\sum_i a(x_i)t_i,
\qquad a\in L^2(X,\nu).
\label{levylin}
\ee
This is exactly the reason of the well-known fact (see, e.g., \cite{Derm90})
that the only L\'evy processes with the so-called
chaotic representation property (which means that the $L^2$
space can be decomposed into the direct sum of the subspaces spanned by ordinary 
multiple stochastic integrals)
are the Gaussian and Poisson processes. 

{\it In the rest of this section we assume that the measure $\Pi$
satisfies the condition (\ref{subord}), i.e., the L\'evy
process is a generalized subordinator, and moreover
the measure $t^2\Pi(t)$ has finite moments of all orders,
and the moment problem for this measure is definite.} In this case we
can obtain a more detailed description of the orthogonal decomposition.

Let $\{P_k(t)\}_{k=0}^\infty$ be the family of orthogonal polynomials on
$\mathbb R_+$ with respect to the measure $t^2\Pi(t)$, and let
$V_{n,k}$
be the subspace spanned by the generalized Charlier
functionals $\frak C^{(n)}_{f_1,\ldots,f_n}$
of order $n$ (see Definition \ref{def:charlier})
corresponding to functions of the form 
$f(x,t)=a(x)tP_{k-1}(t)$, i.e., to functions of the {\it $k$th power} in $t$.

It is not difficult to see that
$$
V_n=\sum_{\la\vdash n}\bigoplus_{(i_1,\ldots,i_k)\atop\mbox{\scriptsize distinct }}
\bigotimes_{j=1}^k V_{\la_j,i_j},
$$
where $\la=(\la_1\ge\ldots\ge\la_k>0)$ is a partition of $n$, and the
tensor product is symmetric.
Thus, rearranging the summands, we obtain the following orthogonal
decomposition for the general L\'evy process:
\begin{multline*}
L^2(D,P_\La)=\bigoplus_n\bigoplus_{\lambda\vdash n}\otimes_sV_{n_k,k}\\
=\mathbb C\oplus V_{1,1}\oplus\left( V_{1,2}\oplus V_{2,1}\right)
\oplus\left (V_{3,1}\oplus (V_{1,1}\otimes_s V_{1,2})\oplus V_{1,3}\right)
\oplus\ldots, 
\end{multline*}
where $\lambda=1^{n_1}2^{n_2}\ldots$
is a partition of $n$ with $n_k$ parts equal to $k$, 
and $\otimes_s$ is the symmetric tensor product. 

The same decomposition can be described in another way
(cf.~\cite{NuSch00}). Given the L\'evy process $\eta=\sum t_i\de_{x_i}$,
consider the processes $\eta_k=\sum_i t_iP_{k-1}(t_i)\de_{x_i}$ (in particular,
$\eta_1=\eta$).
Then $V_{n,k}$ is the space of $n$-multiple stochastic integrals of the process
$\eta_k$.

The corresponding decomposition for the $L^2$ space over the vector-valued
Gaussian process is obtained in a similar way. Namely,
\begin{multline*}
L^2(\widehat{H^\La},\mu^\La)=\bigoplus_n\bigoplus_{\lambda\vdash n}\otimes_s{\cal H}_{n_k,k}\\
=\mathbb C\oplus {\cal H}_{1,1}\oplus\left({\cal H}_{1,2}\oplus {\cal H}_{2,1}\right)
\oplus\left ({\cal H}_{3,1}\oplus ({\cal H}_{1,1}\otimes_s {\cal H}_{1,2})\oplus {\cal H}_{1,3}\right)
\oplus\ldots,
\end{multline*}
where ${\cal H}_{n,k}$ is the subspace spanned by the generalized 
Hermite functionals of order $n$ (see Definition \ref{def:hermite}) 
corresponding to functions of the form
$f(x,t)=a(x)tP_{k-1}(t)$.

\begin{corollary}
In terms of orthogonal decompositions, the canonical isomorphism 
(\ref{isomlevy}) takes the form
$$
\Phi \frak H^{(n)}_{f_1,\ldots,f_n}=
\frak C^{(n)}_{f_1,\ldots,f_n},
$$
where $f_k(x,t)=a_k(x)tP_{k-1}(t)$.
\end{corollary}

\noindent{\bf Example. Gamma processes.}
The standard gamma process on the space
$(X,\nu)$ is the generalized subordinator 
$\ga$ with the L\'evy measure
\be
\La_\Ga(t)=\frac{e^{-t}}{t}dt,\qquad t>0.
\ee
Thus the Laplace transform of the gamma process is given by
\begin{equation}
\Bbb E e^{-\langle a,\ga\rangle}=
\exp\left(-\int_X\log\left(1+a(x)\right)d\nu(x)\right),
\label{galapl}
\end{equation}
where $a$ is an arbitrary nonnegative measurable function on 
$X$ such that
$\int_X\log(a(x)+1)d\nu(x)<\infty$.

To the space $L^2(\ga)$, we can apply all considerations of this section;
observe that the orthogonal polynomials with respect to
$t^2\La_\Ga(t)=te^{-t}dt$ are the Laguerre polynomials with parameter one:
$P_n(t)=L^{(1)}_n(t)$. Note that the Laguerre polynomials appear
in this example not because of the well-known fact that they are orthogonal
with respect to the gamma distribution, which is the infinitely
divisible distribution corresponding to the gamma process, but since
they are orthogonal
with respect to $t^2\La_\Ga(t)$, i.e., to the L\'evy measure
of the gamma process with the density $t^2$.

\section{Representations of the current group $\SL(2,\mathbb R)^X$}
\label{sect:repr}

As an application of the obtained results, we consider the isomorphism
between the Fock space and 
the $L^2$ space over the gamma process (and the isomorphic 
$L^2$ space over the ``infinite-dimensional Lebesgue measure'')
and apply this construction to representations of the current groups
over $\SL(2,\mathbb R)$.

\subsection{The canonical state on $\SL(2,\mathbb R)^X$}

Let $(X,\nu)$ be a standard Borel space with a fixed finite measure $\nu$.
The {\it current group} $G^X=\SL(2, \Bbb R)^X$ on $(X,\nu)$
is the group of Borel bounded $\SL(2,\Bbb R)$-valued
functions on $X$. In other words, $G^X$ consists
of $2\times2$-matrices whose elements are bounded measurable real
functions on $X$. 

The {\it canonical representation} of the current group $G^X$ is
a unitary irreducible representation with spherical function
given by the formula
\begin{equation}
\Om(g(\cdot))=C\exp\left(-\int_X\log
\big(2+\operatorname{Tr}(g(x)g^*(x))\big)
d\nu(x)\right), \qquad g(\cdot)\in G^X.
\label{spherfunc}
\end{equation}
The restriction of this spherical function
to the subgroup of constant
functions (isomorphic to $G=\SL(2,\Bbb R)$) equals
$\Om_0(g)=\frac{C}{2+\operatorname{Tr}gg^*}$,
the so-called {\it canonical state} of $\SL(2,\Bbb R)$, see \cite{GGV73}.

Consider an infinitely divisible positive definite function on a group,
in other words, a continuous one-parameter semigroup of positive
definite functions. The only interesting case is when the generator of
this semigroup is not positive definite, but only 
conditionally positive definite. Then this generator, as
a function on the group, is not bounded, and it is the norm of a nontrivial
cocycle of the group with values in the space of some irreducible 
representation of this group, 
see \cite{GGV73, VerKarp, Gui, Ar70, Par, Shalom}, and others.
The existence of such cocycle is possible only if the identity representation
is not isolated in the space of all irreducible representations
(i.e., if the group does not satisfy Kazhdan's property (T) \cite{Kazhdan}).
Among classical groups, only $SO(n,1)$ and $SU(n,1)$, $n=1,2,\ldots$,
do have this property, and the corresponding cocycle 
and state were found in
\cite{GGV73, GGV74}. It is this state that is called canonical. It 
allows one to define a representation of the current group in the Fock space.
Formula (\ref{spherfunc}) above determines the positive definite function 
on the current group $\SL(2,\mathbb R)^X$ generated by the canonical state;
it is the spherical function generated by the vacuum vector of the corresponding
representation realized in the Fock space. Below we give another
realization of this representation (see Sect.~\ref{sect:lebesgue}).

Note that restrictions of the canonical state to different subgroups
(or commutative subalgebras of the group algebra) determine
different infinitely divisible measures on the dual subgroup
(respectively, the dual space to the algebra), thus diagonalization
of different subgroups or subalgebras generates different infinitely divisible
measures, L\'evy processes, and hence models of the Fock space.

\subsection{The Fock model of the canonical representation}
\label{sect:fock}

Let us describe the Fock model of the canonical representation of 
$\SL(2,\Bbb R)^X$. Consider the so-called {\it special} representation
of the group $G=\SL(2,\Bbb R)$, which is realized in the 
Hilbert space $H=L^2(\Bbb R,\frac{dt}{|t|})$ and is
given by the following formulas:
\begin{eqnarray}
\left(T\left(
 \begin{array}{cc}
 1& 0\\
 b&1
 \end{array}
 \right)\phi\right)(t)&=&e^{ibt}\phi(t),\nonumber\\
 \left(T\left(
  \begin{array}{cc}
  a^{-1}& 0\\
  0&a
  \end{array}
  \right)\phi\right)(t)&=&\phi(a^2t),\\
  \left(T\left(
   \begin{array}{cc}
   0& 1\\
   -1&0
   \end{array}
   \right)\phi\right)(t)&=&\int_{-\infty}^\infty K_0(t,s)\phi(s)ds,\nonumber
   \end{eqnarray}
   where
\begin{equation}
\label{K0}
K_0(t,s)=\frac{1}{2\pi}\int_{-\infty}^\infty
   \frac{1}{|u|^2}e^{-i(tu+su^{-1})}du.
\end{equation}
This is an irreducible unitary representation of discrete series. 
(Note that considered over $\Bbb C$  this representation
is reducible: it decomposes into the sum of two 
irreducible subrepresentations.) However, this representation is distinguished 
as the only representation having a nontrivial cocycle.
\smallskip

Let $\phi_0(t)=e^{-|t|}$ and
fix a cocycle $\beta:G\times H\to H$ given by 
\begin{equation}
\beta(g,t)=T_g\phi_0(t)-\phi_0(t).
\end{equation}

Consider the Hilbert space
$$
H^X=L^2\left(X\times\Bbb R,\, \nu\times\frac{dt}{|t|}\right)
$$
and the corresponding Fock space $\EXP H^X$.
The realization of the canonical representation of 
$G^X=\SL(2,\Bbb R)^X$ in $\EXP H^X$ is given by the
formula
\begin{equation}
U_{g(\cdot)}\EXP{h(\cdot,\cdot)}=
\la(g(x),h(x,t))\cdot\EXP({T_{g(x)}h(x,t)+\beta(g(x),t)}),
\label{fockrepr}
\end{equation}
where
\begin{multline}
\la(g,h)=\exp\left(-\frac12\|\beta\|^2-\langle T_gh,\beta\rangle
\right)\\
=\exp\left(
-\frac12\int_X\int_{\Bbb R}\frac{|\beta(g(x),t)|^2}{|t|}dt\,d\nu(x)
-\int_X\int_{\Bbb R}\frac{T_{g(x)}h(x,t)
\cdot\bar\beta(g(x),t)}{|t|}
dt\,d\nu(x)
\right).\nonumber
\end{multline}
The vacuum vector in the Fock realization is $\EXP 0$, and
the corresponding spherical function equals (\ref{spherfunc}).

Using the isomorphism (\ref{fockgauss}) between the Fock space and the $L^2$
space over the Gaussian white noise, one can
obtain the Gaussian realization of the
canonical representation. We do not reproduce here the 
corresponding formulas, which can be found in \cite{GGV85}.

\subsection{The Lebesgue model of the canonical representation}
\label{sect:lebesgue}

The commutative model of the canonical representation of $G^X$
with respect to the unipotent subgroup was given in \cite{GGV85}.
Another realization of this model, in the $L^2$ space over the
so-called infinite-dimensional Lebesgue measure,
was constructed in \cite{TVY01}. Let us describe this model.

The {\it Lebesgue measure} ${\cal L}^+$ 
on the space $D^+(X,\nu)$
is a $\si$-finite measure
equivalent to the law $\cal G$ of the 
gamma process (see the example at the end of \S\ref{sect:genisom})
with the density given by 
\begin{equation}
\frac{d{\cal L}^+}{d{\cal G}}(\eta)=\exp(\eta(X)).
\label{sigmadens}
\end{equation}
It follows from (\ref{galapl}) and~(\ref{sigmadens}) that the
Laplace transform
of ${\cal L}^+$ equals
\begin{equation}
\int_{D^+}\left[\exp\left(-\int_Xa(x)d\eta(x)\right)\right]d{\cal L}^+(\eta)=
\exp\left(-\int_X\log a(x)d\nu(x)\right).
\label{sigmalaplace}
\end{equation}
The {\it Lebesgue measure
on $D(X)$} is the convolution ${\cal L}^+*{\cal L}^-$, where
${\cal L}^-$ is the image of
${\cal L}^+$ under the mapping $\eta\to-\eta$.

An arbitrary measurable function $a:X\to\mathbb R_+$ with
$\int_X|\log a(x)|d\nu(x)<\infty$ defines a multiplicator
$M_a:D\to D$ by the formula
$$
M_a:\eta=\sum_it_i\de_{x_i}\mapsto \sum_i a(x_i)t_i\de_{x_i}.
$$
As shown in \cite{TVY01}, the Lebesgue measure is projective invariant 
with respect to the group of multiplicators, namely,
$$
\frac{dM_a({\cal L})}{d{\cal L}}=
\exp\left(-\int_X\log a(x)d\nu(x)\right).
$$
This key property of the infinite-dimensional Lebesgue measure is
a consequence of a remarkable quasi-invariance property of the gamma process
(see \cite{TV99, TVY01}). In particular, it makes it possible to 
construct a representation of the current group in the
$L^2$ space over the Lebesgue measure. Note also that though there exist
other subordinators quasi-invariant with respect to the group of
multiplicators (see \cite{LiSh}), however, the gamma process is the only
subordinator that admits an equivalent measure that is projective
invariant with respect to this group
(\cite{TVY01}).

Consider the triangular subgroup $\cal T$ of $\SL(2,\mathbb R)^X$:
$$
{\cal T}=\left\{
T_{a,b}=\left(
 \begin{array}{cc}
a(\cdot)^{-1}& 0\\
b(\cdot)&a(\cdot)
\end{array}
\right)\right\}.
$$

\begin{theorem}[\cite{TVY01}]
The formula
\begin{equation}
\label{lebesguerepr}
{\cal U}(T_{a,b})F(\eta)= 
\exp\left(\int_X\log |a(x)|d\nu(x)+i\int_X a(x)b(x)d\eta(x)\right)F(M_{a^2}\eta)
\end{equation}
defines a unitary irreducible representation of the triangular
subgroup $\cal T$ in the space $L^2(D,{\cal L})$,
which is extendable to a unitary irreducible representation of
the whole group $\SL(2,\mathbb R)^X$.
\label{th:repr}
\end{theorem}

\subsection{Isomorphism of the Fock and Lebesgue models of the
canonical representation}
\label{sect:isomrepr}

The Fock model (\ref{fockrepr}) and the Lebesgue model (\ref{lebesguerepr}) define
isomorphic representations, since their spherical functions 
coincide. However, now we can use the canonical isomorphism between the space
of square integrable functionals of the gamma process and the Fock space to
construct explicitly the isomorphism
of these realizations.

\begin{theorem}\label{th:isomrepr}
The isometry of the spaces $\EXP H^X$ and $L^2(D,{\cal L})$ that
intertwines the Fock realization $U$ and the Lebesgue realization $\cal U$ 
of the canonical representation of the current group
$\SL(2,\mathbb R)^X$
is given by
$$
\EXP h\leftrightarrow\Psi_h,\qquad 
h\in L^2(X\times\mathbb R,\nu\times\frac{dt}{|t|}),
$$
where
\be
\label{45}
\Psi_h(\eta)=\prod_i\left(h(x_i,t_i)+e^{-|t_i|/2}\right)\cdot
\exp\left(-\int_X\int_{\Bbb R}\frac{h(x,t)\cdot e^{-|t|/2}}{|t|}dt\,d\nu(x)
\right)
\ee
for $\eta=\sum t_i\de_{x_i}\in D$. 
\end{theorem}
\begin{proof}
Formula (\ref{45}) defines an isometry between $\EXP H^X$ and $L^2(D,{\cal L})$,
as follows from Theorem~\ref{th:levy}
in the special case of the gamma process and the obvious isometry
between $L^2(D,\cal G)$ and $L^2(D,{\cal L})$ given by
$F(\eta)\leftrightarrow F(\eta)e^{-\eta(X)/2}$. It is not difficult to 
verify by direct calculation that this isomorphism intertwines the 
representations $U$ and $\cal U$.
\end{proof}

Note that the vacuum vector
in the Lebesgue realization is
$\Psi_0(\eta)=e^{-\frac{|\eta|(X)}{2}}$, where $|\eta|=\sum |t_i|$
is the total charge of the (signed) measure $\eta$.

\begin{corollary}
The action of the involution $\si=\left(
\begin{array}{cc}
0&1\\
-1&0
\end{array}
\right)$ in the Lebesgue realization of the canonical representation
of $G^X$ is given by the formula
\begin{equation}
U_\si\Psi_{f(\cdot,\cdot)}=\Psi_{T_\si f(\cdot,\cdot)},
\end{equation}
where 
\begin{equation}
T_\si f(x,t)=\int_{-\infty}^\infty K_0(t,s)f(s)ds,
\end{equation}
with the kernel $K_0$ given by (\ref{K0}).
\end{corollary}
\begin{proof} Follows from Theorem~\ref{th:isomrepr} and (\ref{fockrepr}), since
$\beta(\si,t)\equiv0$.
\end{proof}

\section*{Appendix}
\addcontentsline{toc}{section}{Appendix}

\subsection*{A. An example of a zero-dimensional non-Fock factorization
(a model of hierarchical voting \cite{TsirVer98})}
\addcontentsline{toc}{subsection}{A. An example of a zero-dimensional non-Fock factorization}
In this appendix, we reproduce the example of a non-Fock zero-dimensional
factorization constructed in
\cite{TsirVer98}.

We consider the simplest, in fact purely combinatorial, model of a 
Hilbert and measure factorization, over a Cantor compactum, that is
not isomorphic to a Fock factorization (i.e., is not linearizable).
This model is determined by two positive integers
$m,r>1$ and a symmetric map $\phi:X_r^m \to X_r$,
where $X_r$ is a set consisting of $r$ elements (it is convenient to enumerate them
by the numbers $0,1, \dots, r-1$) and $X_r^m$ is its  $m$th power, and
the number of points in the preimage of each point 
$x \in X_r$ is the same, i.e., $\#(\phi^{-1}(x))=m^{r-1}$.
The latter condition implies that the 
$\phi$-image of the uniform measure on
$X_r^m$ is the uniform measure on $X_r$.

For each such triple  $(m,r,\phi)$, we will construct a factorization;
under very wide assumptions
on $\phi$, these factorizations are not isomorphic to a Fock factorization
and have a large group of symmetries.

The map $\phi$ is called {\it antiadditive} (respectively, 
{\it antimultiplicative}) if it satisfies the following condition.
If for a function $g:X_r \to \mathbb C$, there exists a function
$f:X_r \to \mathbb C$ such that the following relation holds
identically (i.e., for any 
$a_1 \in X_r, \dots, a_m\in X_r$):
\be\label{*}
f(\phi (a_1, \dots a_m))=g(a_1)+ \dots + g(a_m)
\ee
(respectively,
\be\label{**}
f(\phi (a_1, \dots a_m))=g(a_1)\cdot \dots \cdot g(a_m)),
\ee
then the function $f$ (and hence $g$) is a constant.

For example, if $X_r$ is an additive or multiplicative group, and
$\phi$ is the group operation, then nonconstant solutions of these
equations are additive or multiplicative characters. Here are examples
of antiadditive maps.

\smallskip\noindent{\bf Examples. 1.}
The model of voting by majority:
$m=3$, $r=2$, and $\phi$ is given by
$$
\phi (a,b,b)=b, \qquad a,b =0,1. 
$$

\noindent{\bf 2.} Let $m=2$, $r=3$ (i.e., $X_r=\{0,1,2\}$), 
and let $\phi$ be given by the table of values
$$
\phi=\begin{array}{ccc}
  2&2&0\\
  2&0&1\\
  0&1&1
\end{array}
$$

It is easy to see that in both examples 
(\ref{*}) and (\ref{**}) have no nonconstant solutions.

The paper \cite{TsirVer98} contains a convenient criterion for 
the solutions of (\ref{*}) and (\ref{**}) to be constant functions
(see below). It turns out that this case is generic,
only in degenerate cases (similar to group laws) nonconstant
solutions appear. In our construction, 
the absence of such solutions will guarantee
the absence of additive and multiplicative vectors in the constructed 
factorization. The key role in the sequel is played by the
following condition on the map $\phi$.

\smallskip\noindent{\bf Abundance condition.}
Let we are given a map  $\phi: X_r^m \to X_r$. Fix 
$m-1$ arguments in an arbitrary way
(due to the symmetry, it does not make difference what
arguments we choose), take all maps of the set
$X_r$ into itself obtained in this way:
$a_m \mapsto \phi_{a_1, \dots a_{m-1}}(a_m)\equiv \phi(a_1, \dots a_m)$,
and consider the subsemigroup generated by all
these maps in the semigroup of all maps of the set $X_r$ 
into itself.

\begin{definition}
The map
$\phi$ is called {\em abundant} if the obtained subsemigroup contains
at least one constant map.
\end{definition}

It is easy to check that abundance is a generic condition. For example,
for $m=2$, it does not hold only for those maps 
$\phi$ that determine a semigroup law on the set 
$X_r$; in this case formula (\ref{**}) defines a multiplicative character
of the group or semigroup. The abundance condition
also appears in the theory of Markov chains. 

Before constructing a Hilbert factorization for an arbitrary triple
  $(m,r,\phi)$, let us describe the corresponding probability space.

Let  $T_m$ be the infinite rooted $m$-ary tree, and assume that each its
vertex is assigned a random variable that takes $r$ values with equal probabilities,
the random variables of the same level (i.e., at the same distance from
the root) being independent and the random variable
$\xi_v$ corresponding to a vertex $v$ being equal to
$\phi(\xi_{v_1},\dots, \xi_{v_m})$, where $v_1, \dots, v_m$ are the 
sons of $v$, and $\phi$ is the map
(``voting'') defined above.

The probability space
 $\Omega=\Omega(m,r,\phi)$ is the space of realizations of this family 
 of random variables, i.e., the space of all functions
$f$ on the set of vertices of the tree $T_m$ with values in the set
$X_r=\{1,2, \ldots, r \}$ that satisfy the above condition:
$\phi(f(v_1), f(v_2), \ldots, f(v_m)) = f(v)$,
where $v_1, \dots, v_m$
are the sons of the vertex  $v$. (Following our analogy, one may call this space
the space of ballot-papers).\footnote{Note that in the ``model of 
voting by majority'' from example~1 above,
$2^n$ voters of the $n$th level can legally defeat all
$3^n$ voters participating in the vote; thus already for the two-level system
($n=2$, the total number of voters is 9), four voters can
defeat the remaining five voters, though the probability of this event is small.}
By the properties of the map $\phi$, this space is equipped with a well-defined
uniform measure, and
the values of functions at different vertices of the same level
(the voters of the same level) are independent with respect to this measure.
Note that the space
$\Omega$ with the uniform measure is the inverse limit of the finite spaces
$X_r^m$ with the uniform measures with respect to the projections defined
by $\phi$.

It is useful to give another interpretation of the space
$\Omega$. Let $K_m$ be the Cantor compactum of all infinite paths in the tree
$T_m$ endowed with the natural totally disconnected topology. Each function
$f\in \Omega$ determines a {\it pseudomeasure} $\nu_f$
on the cylinder sets of the space $K_m$. Namely, 
by definition, the value of the pseudomeasure
$\nu_f$ on the cylinder $C_v$ of all paths going through a vertex $v$
is equal to  the value of the function $f$ at the vertex $v$. 
By a pseudomeasure, we mean a function $\nu$ defined on the algebra
of cylinder sets of the compactum $K_m$ and satisfying a unique condition
on the values at elementary cylinders\footnote{An elementary cylinder
of order $k$ is the set of paths in the tree with a given initial
segment of length $k$.}, which reproduces the condition on the functions of
the space $\Om$:
$\nu(C_v)=\phi(\nu(C_{v_1}),\ldots, \nu(C_{v_m}))$; one may call this condition
$\phi$-additivity. Thus we have described the space
$\Omega$ as a space of pseudomeasures. The measure on this space allows
us to speak about random pseudomeasures.

The space $\Omega$ has a natural measure factorization in the sense of  
Definition~\ref{def:measfact} over the Boolean algebra of cylinder sets
in the space $K_m$. Note that the map
$\zeta$ from the definition of a measure factorization is defined in our case only
on elementary cylinders, however, it can be correctly extended to the 
Boolean algebra (but not the $\si$-algebra!) generated by cylinders,
since every cylinder can be uniquely decomposed into 
elementary ones. However, the $\phi$-additivity condition must
hold only for the decomposition of an elementary cylinder into elementary ones.

Now we are in a position to define a Hilbert factorization, which 
will be non-Fock under a certain condition on the map $\phi$.
But first let us give the following analogy. The above description is
similar to the following  nonconventional simple description of
processes with independent values, namely, the approximative description.
For simplicity, we will speak only about the white noise and use the 
same notation as before. In our example, replace a finite space
$X_r$ by the real line
$\mathbb R$ with the standard Gaussian measure, and let  
$\phi:\mathbb R^m \to\mathbb R$ be the normalized sum:
  $\phi(\xi_1, \dots, \xi_m)=\frac{\xi_1 +\ldots  + \xi_m}{\sqrt m}$ 
($\phi$ sends the standard Gaussian measure on $\mathbb R^m$ 
to the Gaussian measure on $\mathbb R$).
It is not difficult to see that in this case our construction
leads to a space $\Omega$ whose elements are ordinary additive measures
on the same Cantor compactum $K_m$, and the Gaussian probability measure on
$\Omega$ is defined by the condition that the value of a (random)
additive real-valued measure on every cylinder is the integral 
of a realization of the standard white noise with the base space
$K_m$ over this cylinder. 
In other words, we have represented the Gaussian measure determined by the
white noise as the inverse limit (in the sense of linear spaces)
of the Gaussian measures on $\mathbb R^n$.
 
 In some sense, in our example, the value of a pseudomeasure
on a cylinder of the set $K_m$ can be also regarded as the result of measuring
a certain nonlinear noise (``black noise'',
though this term looks too gloomy) on this cylinder; but instead of
additivity we have only the 
``$\phi$-additivity'' defined above.

Now it is not difficult to describe the Hilbert space
$L^2(\Omega)$ and explain the appearance of a non-Fock factorization.
The Hilbert space $H=L^2(\Omega)$ is of course the inductive limit in
$k$ of the finite-dimensional spaces
$H_k=L^2(X_r^{m^k})$ with respect to the embeddings determined by the map
$\phi$; here $H_k$ is the space of functions on the product of the spaces
$X_r$ over the vertices of the $k$th level. But the embedding
$$
\bar \phi_k:  H_k  \to H_{k+1}
$$ 
is the tensor product of $k$ copies of the embedding
$\bar\phi_1: L^2(X_r) \to L^2(X_r^m)$, thus it suffices to define
only the latter embedding; it is defined on the basis 
$e_1, \ldots,e_r$ of the space $L^2(X_r)$ by the formula
 $\bar \phi_1(e_i)=\chi_{\phi^{-1}(i)}$, where $\chi_E$
 is the characteristic function of the set
$E$, and $\phi^{-1}(i)$ is the preimage of an element
 $i \in X_r$ under the map
 $\phi$. The constructed embeddings are obviously isometric,
and they define a scalar product in the inductive limit of spaces.
Recall that the inductive limit in the category of Hilbert spaces
is the completion of the algebraic limit with respect to the Hilbert norm.
It follows immediately from the above considerations that the constructed
inductive limit can be naturally identified with the space
$L^2(\Omega)$. It also follows from construction that the following
proposition holds. 

\begin{proposition}
The space $L^2(\Omega)$ has a factorization over the Boolean
algebra of cylinder sets of the space $K_m$.
\end{proposition}

Recall that the continuity of
a factorization at a point means the following: for each
decreasing sequence of cylinders whose intersection
is a point, the corresponding 
sequence of operator algebras converges to the algebra consisting of
scalar operators.

\begin{theorem}
The constructed factorization is continuous at a point if and only if
the map $\phi$ is abundant. Under this condition, the factorization is
non-Fock.
\end{theorem}

\begin{proof}
Indeed, the abundance condition implies the antiadditivity, i.e., the
absence of additive, and hence (by Theorem \ref{th:log}) 
multiplicative vectors.
\end{proof}

\begin{corollary}
The constructed factorization cannot be extended to the Boolean algebra
of $\bmod 0$ classes of measurable subsets of the space
$K_m$ with the natural (Lebesgue) product measure $\mu$.
\end{corollary}

Indeed, by Tsirelson's theorem, such extension must be a Fock
factorization. On the other hand, a direct calculation shows
that not for every convergent sequence of cylinder sets, 
the limit of the corresponding random variables does exist.

It follows from the definition of the constructed factorization
that its group of symmetries includes the whole group 
of automorphisms of the tree and, moreover, contains also other transformations,
see \cite{TsirVer98}.  The study of properties
of these factorization, in particular, the structure of the Hilbert
space equipped with such factorization is of great interest. 
It is very important to
consider operators similar to the canonical operators
in a Fock factorization. This space may be eventually useful
for the representation theory of groups and
$C^*$-algebras. From the probabilistic viewpoint, the example
under consideration is apparently related to the theory of
branching processes.

The constructed factorization has a zero-dimensional base
(the Cantor compactum). The construction of a non-Fock factorization with
a one-dimensional base is more complicated, see
\cite{TsirVer98}. The difficulty appears since we must
coordinate the factorization at different representations of the interval
as the union of intervals (while the decomposition of elementary 
cylinders of the Cantor compactum is unique). However, in the one-dimensional
case, the factorization is also constructed by means of the inverse limit
of Gaussian measures with nonlinear, as above, projections.
Examples of non-Fock factorizations in dimensions greater
than one are not yet constructed.

\subsection*{B. On the spectrum of a Fock factorization}
\addcontentsline{toc}{subsection}{B. On the spectrum of a Fock factorization}

Consider an arbitrary factorization over a Boolean algebra of sets.
Associate with each measurable set $B$ of the ``base'' the projection
(conditional expectation) $P_B$ onto the subspace of functionals 
that depend only on the restriction of a realization of the process to this set. Obviously, these
projections commute with each other. Thus they generate a commutative
$C^*$-algebra
$\frak M$. The following notion is introduced by Tsirelson
\cite{Tsir98}.

\begin{definition}
The spectrum of the commutative $C^*$-algebra $\frak M$ 
is called the {\em spectrum of the factorization}.
\end{definition}

By definition, the spectrum is an invariant of a factorization. Consider
the spectrum of the Fock factorization generated by the one-dimensional
Gaussian process. Let us prove that in this case the algebra
$\frak M$ is maximal in the algebra of all operators, i.e.,
each operator that commutes with all these projections belongs
to the same algebra. First of all, the projections to the ``chaos''
of a given order commute with all
$P_B$, thus we may restrict the subalgebra
$\frak M$ to the subspace of the given chaos, and this restriction 
is a maximal commutative subalgebra in the subalgebra of operators
of this chaos. Finally, it is obvious that two different chaoses are
separated by the algebra
$\frak M$. 

\begin{theorem}
The spectrum of the one-dimensional Fock factorization is the disconnected
union of symmetrized $n$-tuples from the base equipped with the
symmetrized powers of the measure on the base.
\end{theorem}

For a multidimensional Fock factorization, the algebra
$\frak M$ is not maximal, but the spectrum is the same.

\subsection*{C. Cycle index}
\addcontentsline{toc}{subsection}{C. Cycle index}
Let $\sn$ be the symmetric group of degree $n$. Given a permutation $g\in\sn$,
denote by  $c_k(g)$ the
number of cycles of length $k$ in $g$. Let
$t=(t_1,t_2,\ldots)$ be a sequence of indeterminates.
The (augmented) cycle index of the symmetric group $\sn$ is
\be
\label{cyndex}
\ti Z(\sn)=\ti Z(\sn)[t]=\sum_{g\in\sn}t_1^{c_1(g)}t_2^{c_2(g)}\cdot\ldots.
\ee
A well-known formula (see, e.g., \cite[(5.30)]{St}) claims that
\be
\label{cycleindex}
\sum_{n=0}^\infty \tilde Z(\sn)[t]\frac{z^n}{n!}=
\exp\sum_{i=1}^\infty t_i\frac{z^i}{i}.
\ee

\subsection*{D. Orthogonal polynomials}
\addcontentsline{toc}{subsection}{D. Orthogonal polynomials}

For convenience, in this Appendix we reproduce
necessary formulas concerning the classical orthogonal polynomials
of Hermite and Charlier, which play an important role in
the theory of orthogonal decompositions for L\'evy processes.
A standard reference on the theory of orthogonal polynomials is
the monograph \cite{Szego62}. For combinatorial aspects of orthogonal
polynomials, see also \cite{St}.

\subsubsection*{D.1. Hermite polynomials}

\noindent{\bf Definition:}
\begin{eqnarray}
H_n(x)&=&H_n^1(x)=(-1)^ne^{\frac{x^2}{2}}\frac{d^n}{dx^n}e^{-\frac{x^2}{2}}
\nonumber,\\
H_n^a(x)&=&a^{\frac n2}H_n(x/\sqrt a).
\label{hermdef}
\end{eqnarray}

\noindent{\bf Orthogonality:}
orthogonal on $\mathbb R$ with respect to the normal law $N(0,a)$ with zero mean
and variance $a$:
\be
\frac{1}{\sqrt{2\pi a}}\int_{-\infty}^\infty 
H_n^a(x)H_m^a(x)e^{-\frac{x^2}{2a}}\,dx=
\de_{nm}a^nn!.
\ee

\noindent{\bf Generating function:}
\be
\label{hermgen}
\sum_{n=0}^\infty \frac{H_n^a(x)t^n}{n!}=e^{tx-\frac{t^2a}{2}}.
\ee

\noindent{\bf Combinatorial description:}
\be
H_n^a(x)=\sum_{g\in\inv_n} x^{c_1(g)}(-a)^{c_2(g)}=
\ti Z_{\sn}[x,-a,0,0,{\ldots}],
\label{hermcomb}
\ee
where $\inv_n$ is the set of all involutions in $\sn$.

\subsubsection*{D.2. Charlier polynomials}

\noindent{\bf Definition:}
\be
C^a_n(x)=a^n\sum_{j=0}^n(-1)^{n-j}{n\choose j}a^{-j}{x\choose j}j!.
\label{chardef}
\ee

\smallskip
\noindent{\bf Orthogonality:}
orthogonal on $\mathbb Z_+$ with respect to the Poisson law $P_a$
with parameter $a$:
\be
e^{-a}\sum_{k=0}^\infty C^a_n(k)C^a_m(k)\frac{a^k}{k!}=\de_{nm}a^nn!.
\ee

\noindent{\bf Generating function:}
\be
\label{chargen}
\sum_{n=0}^\infty\frac{C^a_n(y)t^n}{n!}=(1+t)^ye^{-ta}.
\ee

\noindent{\bf Combinatorial description:}
\be
C_n^a(x)=\sum_{g\in\sn} (x-a)^{c_1(g)}x^{-c_2(g)+c_3(g)-c_4(g)+\ldots}
=\ti Z_{\sn}[x-a,-x,x,-x,{\ldots}].
\label{charcomb}
\ee

The Charlier polynomials satisfy the following convenient formula:
\be
\label{chform}
\frac{(-1)^n}{a^n}C^a_n(k)=\frac{(-1)^k}{a^k}C^a_k(n).
\ee

\section*{Notation}
\addcontentsline{toc}{section}{Notation}
\begin{tabular}{lp{20cm}}
$(X,\nu)$ & standard Borel space with a continuous finite measure \\
$\al$ & standard Gaussian white noise on $(X,\nu)$ \\
$\pi$ & Poisson process on $(X,\nu)$\\
$\Pi$ & L\'evy--Khintchin measure of a L\'evy process \\
$\Phi$ & canonical Poisson--Gauss isomorphism\\
$\frak H^{(n)}_{f_1,\ldots,f_n}$ & generalized Hermite functional\\
$\frak C^{(n)}_{f_1,\ldots,f_n}$ & generalized Charlier functional\\
$D$ & space of finite real discrete measures on $X$\\
$\frak S_n$ & symmetric group of degree $n$\\
$\inv_n$ & the set of all involutions in $\frak S_n$
\end{tabular}

%\bibliographystyle{unsrt}
%\bibliography{engisom}

\end{document}